\documentclass[twocolumn]{autart}    

\usepackage{enumitem}
\usepackage{mathtools} 
\usepackage{dsfont, url}
\usepackage{newtxtext} 
\usepackage{textcomp} 
\usepackage[varg,bigdelims]{newtxmath}
\usepackage{MnSymbol}
\usepackage[usenames, dvipsnames]{xcolor}
\usepackage[%
	colorlinks={true},%
	citecolor={blue},%
	urlcolor={blue},%
	linkcolor={blue}%
]{hyperref}
\def\revised#1{{\color{black}#1}}
\definecolor{mycolor1}{RGB}{0,102,102}%
\definecolor{mycolor2}{RGB}{181,179,246}%
\definecolor{mycolor3}{RGB}{153,255,153}%
\definecolor{mycolor4}{rgb}{0,0.55,0.55}%
\usepackage{cite}
\DeclareMathOperator*{\minimize}{minimize}
\DeclareMathOperator*{\sbjto}{subject\ to}
\DeclareMathOperator*{\blkdiag}{blkdiag}

\DeclareMathOperator{\rank}{rank}
\DeclareMathOperator{\sat}{sat}
\DeclareMathOperator{\trace}{tr}

\renewcommand{\le}{\leqslant}
\renewcommand{\leq}{\leqslant}

\renewcommand{\geq}{\geqslant}

\newcommand{\R}{\mathds{R}}
\newcommand{\Nz}{\mathds{N}}

\newcommand{\EE}{\mathds{E}}

\newcommand{\inprod}[2]{\left\langle{#1}, {#2}\right\rangle}
\newcommand{\pmat}[1]{\begin{bmatrix}#1\end{bmatrix}}
\newcommand{\abs}[1]{\left|#1\right|}
\newcommand{\norm}[1]{\left\|#1\right\|}
\newcommand{\secref}[1]{\S\ref{#1}}

\newcommand{\transp}{^\top}

\newcommand{\ones}{\mathds{1}}
\newcommand{\zeros}{\textbf{0}}
\newcommand{\st}{x}
\newcommand{\stortho}{\st^o}
\newcommand{\stinit}{\overline{\st}}

\newcommand{\A}{A}
\newcommand{\calA}{\mathcal{A}}
\newcommand{\Aortho}{\A_o}
\newcommand{\Aschur}{\A_s}
\newcommand{\B}{B}
\newcommand{\calB}{\mathcal{B}}
\newcommand{\calD}{\mathcal{D}}
\newcommand{\Bortho}{\B_o}
\newcommand{\Bschur}{\B_s}

\newcommand{\control}{u}
\newcommand{\controlset}{\mathds{U}}
\newcommand{\wnoise}{w}

\newcommand{\cnoise}{\nu}

\newcommand{\calQ}{\mathcal{Q}}
\newcommand{\calR}{\mathcal{R}}

\newcommand{\calS}{\mathcal{S}}
\newcommand{\calG}{\mathcal{G}}
\newcommand{\calL}{\mathcal{L}}
\newcommand{\calM}{\mathcal{M}}

\newcommand{\costps}{c_{\mathrm{s}}}
\newcommand{\costfinal}{c_{\mathrm{f}}}

\newcommand{\regular}{\varrho}

\newcommand{\lra}{\longrightarrow}

\newcommand{\ee}{\mathfrak{e}}
\newcommand{\offset}{\boldsymbol{\eta}}
\newcommand{\weight}{\boldsymbol{\Lambda}}
\newcommand{\gain}{\boldsymbol{\Theta}}

\newcommand{\reachindex}{\kappa}
\newcommand{\reachab}{\mathrm{R}}
\newcommand{\Let}{\coloneqq}
\newcommand{\teL}{\eqqcolon}

\newtheorem{assumption}{Assumption}

\newtheorem{remark}{Remark}

\usepackage{pgf,tikz}
\usepackage{adjustbox}
\usepackage{verbatim} 
\usepackage{xargs}
\usetikzlibrary{plotmarks}
\usetikzlibrary{chains}
\usepackage{pgfplots}
\pgfplotsset{compat=newest}
\pgfplotsset{plot coordinates/math parser=false}
\newlength\figureheight
\newlength\figurewidth
\usetikzlibrary{shapes,arrows}
\usetikzlibrary{positioning}
\usetikzlibrary{decorations,decorations.pathmorphing,decorations.pathreplacing}
\usetikzlibrary{shadows,trees}
\pgfdeclaredecoration{switch}{initial}
{
    \state{initial}[width=+\pgfdecoratedinputsegmentremainingdistance]
    {
        \pgfsetlinewidth{1pt}
        \pgfpathcircle{\pgfpoint{0}{0}}{0.05\pgfdecoratedinputsegmentremainingdistance}
        \pgfpathlineto{\pgfpoint{0.25\pgfdecoratedinputsegmentremainingdistance}{0}}
        \pgfpathcircle{\pgfpoint{0.25\pgfdecoratedinputsegmentremainingdistance}{0}}%
            {0.05\pgfdecoratedinputsegmentremainingdistance}
        \pgfpathlineto{\pgfpoint{0.75\pgfdecoratedinputsegmentremainingdistance}{\pgfdecorationsegmentlength}}
        \pgfpathcircle{\pgfpoint{0.75\pgfdecoratedinputsegmentremainingdistance}{0}}%
            {0.05\pgfdecoratedinputsegmentremainingdistance}
        \pgfpathlineto{\pgfpoint{\pgfdecoratedinputsegmentremainingdistance}{0}}
        \pgfpathcircle{\pgfpoint{\pgfdecoratedinputsegmentremainingdistance}{0}}%
            {0.05\pgfdecoratedinputsegmentremainingdistance}
    }
    \state{final}
    {
        \pgfpathmoveto{\pgfpointdecoratedpathlast}
    }
}
\allowdisplaybreaks
\begin{document}

\begin{frontmatter}

\title{Sparse and Constrained Stochastic Predictive Control for Networked Systems}

\thanks[footnoteinfo]{D.\  Chatterjee was supported in part by the grant 12IRCCSG005 from IRCC, IIT Bombay.}
\thanks[footnoteinfo]{Authors are thankful to Prof. Masaaki Nagahara and Prof. D. Manjunath for helpful discussions.}
\thanks[footnoteinfo]{ ~ A preliminary version of parts of this work has been presented in $10^{th}$ IFAC symposium on Nonlinear Control Systems, NOLCOS 2016.}
\author[IIT]{Prabhat K. Mishra} 
\author[IIT]{Debasish Chatterjee}
\author[Paderborn]{Daniel E. Quevedo}

\address[IIT]{Systems \& Control Engineering,
Indian Institute of Technology Bombay, 
India.  \tt{prabhat@sc.iitb.ac.in}, \tt{dchatter@iitb.ac.in}}
\address[Paderborn]{Department of Electrical Engineering (EIM-E), Paderborn University
, Germany. \tt{dquevedo}@ieee.org}


\begin{abstract}                          
This article presents a novel class of control policies for networked control of Lyapunov-stable linear systems with bounded inputs. The control channel is assumed to have i.i.d. Bernoulli packet dropouts and the system is assumed to be affected by additive stochastic noise. 
Our proposed class of policies is affine in the past dropouts and saturated values of the past disturbances. We further consider a regularization term in a quadratic performance index to promote sparsity in control. We demonstrate how to augment the underlying optimization problem with a constant negative drift constraint to ensure mean-square boundedness of the closed-loop states, yielding a convex quadratic program to be solved periodically online. The states of the closed-loop plant under the receding horizon implementation of the proposed class of policies are mean square bounded for any positive bound on the control and any non-zero probability of successful transmission. 
\end{abstract}

\end{frontmatter}

\section{Introduction}

\par An ever-increasing number of modern control technologies requires remote computation of control values that are then transmitted to the actuators over a network. Examples include heat, ventilation, and air-conditioning systems (HVAC) \cite{HVAC_survey, kelmanHVAC, hvac2013, building2014} and cloud-aided vehicle control systems \cite{li2014cloud, li2015h, vehicle2015, vehicleColin2015}. In all such systems, a crucial role is played by the transmission channel and the communication protocol employed for the transmission of control commands.   
Due to fading and interference, the transmitted control commands may be delayed or lost, thereby affecting both qualitative and quantitative properties of the system. Since in networked systems rate limited channels are shared among various devices, sparse controls are also desirable and tractability is essential for the implementation. 
Moreover, from an operational stand point, in almost all practical applications there are hard constraints on the controls, and standard control design methods do not apply directly. Furthermore, since stability is one of the most desirable features, it is important to guarantee stability in the context of imperfect communication channel, stochastic noise and bounded controls. This article proposes a sparse, computationally tractable, constrained and stabilizing networked control method for stochastic systems based on predictive control techniques. 

\par
\begin{figure}[t]
\centering
\begin{adjustbox}{width = \columnwidth}
%

\tikzset{
  basic/.style  = {draw, text width=2cm, drop shadow, font=\sffamily, rectangle},  
  final/.style  = {draw, text width=1cm, drop shadow, font=\sffamily, circle},
  root/.style   = {draw, text width=4cm, drop shadow, font=\sffamily, rectangle, rounded corners=2pt, thin, align=center,
                   fill=red!10},
  level 2/.style = {basic, thin,align=center, fill=red!30,
                   text width=2cm},
  level 3/.style = {basic, thin, align=center, fill=green!20, text width=2cm},
  level 4/.style = {basic, thin, align=center, fill=green!40, text width=3cm},
  level 5/.style = {basic, thin, align=center, fill=green!60, text width=5em},
  level 6/.style = {basic, thin, align=center, fill=yellow!20, text width=1.2cm},
  level 7/.style = {final, thin, align=center, fill=blue!10, text width=1cm},
level 8/.style = {draw, text width=2cm, drop shadow, font=\sffamily, rectangle, rounded corners=3pt, thin, align=center,
                   fill=blue!40},
  level 9/.style = {basic, thin, align=center, fill=green!90, text width=1cm}
}

\begin{tikzpicture}[
  level 1/.style={sibling distance=62mm},
  edge from parent/.style={->,draw},
  >=latex]

\node[root] {Stabilizing constrained control over networks}
  child {node[level 2] (tractable) {tractability}}
  child {node[level 2] (sparse) {sparsity}};

\begin{scope}[every node/.style={level 3}]
\node [below of = tractable] (predictive control) {predictive control};

\end{scope}
 \begin{scope}[every node/.style={level 4}]
 \node [right of = predictive control, xshift = 2.3cm] (decision variable) {decision variables};
 \node [below of = predictive control, xshift=-35pt] (cost) {performance index};
 \end{scope}
 
 \begin{scope}[every node/.style={level 5}]
 \node [below of = cost, xshift=-35pt] (expected) {expected};
 \node [below of = cost, xshift=40pt, yshift = -0.2cm] (deterministic) {certainty-equivalent};
 
 \node [below of = decision variable, xshift=-35pt] (sequence) {sequences};
 \node [below of = decision variable, xshift=40pt] (policy) {policies};

 \end{scope}
\begin{scope}[every node/.style={level 6}]
%
\node [below of = policy, xshift=-35pt] (wnoise) {process noise};
 \node [below of = policy, xshift= 45pt] (cnoise) {dropouts};
 \end{scope}
\begin{scope}[every node/.style={level 7}]
 \node [below of = wnoise, yshift = -1 cm] (SMPC) {SPC};
\node [below of = sequence, yshift= -1cm] (PPC) {PPC};
\end{scope}

 \begin{scope}[every node/.style={level 8}]
   \node [below of = cnoise, yshift= -1cm] (New) {proposed approach};
 \end{scope}
 
 
\draw[->] (tractable.south) -- (predictive control.north);
\draw[->] (predictive control.east) -- (decision variable.west);
\draw[->] (predictive control.west) -| (cost.north);
\draw[->] (decision variable.south) -- (sequence.north);
\draw[->] (decision variable.south) -- (policy.north);

\draw[->] (cost.south) -- (expected.north);
\draw[->] (cost.south) -- (deterministic.north);
\draw[->] (expected.south) |- (SMPC.west);
\draw[->] (policy.west) -| (wnoise.north);
\draw[->] (policy.east) -| (cnoise.north);
\draw[->] (wnoise.south) -- (SMPC.north);
\draw[->, thick, red] (cnoise.south) -- (New.north);
\draw[->, thick, red] (SMPC.east) -- (New.west);
\draw[->, thick, red] (sparse.345) -- (New.27);
\draw[->] (sequence.south) -- (PPC.north);
\draw[->] (deterministic.south) |- (PPC.west);


\end{tikzpicture}
\end{adjustbox}
\caption{The approach proposed in the present article extends the results of stochastic predictive control (SPC) as proposed in \cite{ref:PDQ-15} by incorporating communication imperfection models and protocols explicitly into the controller design. In particular, past dropouts are considered in the feedback policy to formulate an sparsity-promoting optimization program. Here, PPC stands for packetized predictive control as described in \cite{quevedo-12}.}
\label{Fig:workflow}
\end{figure}
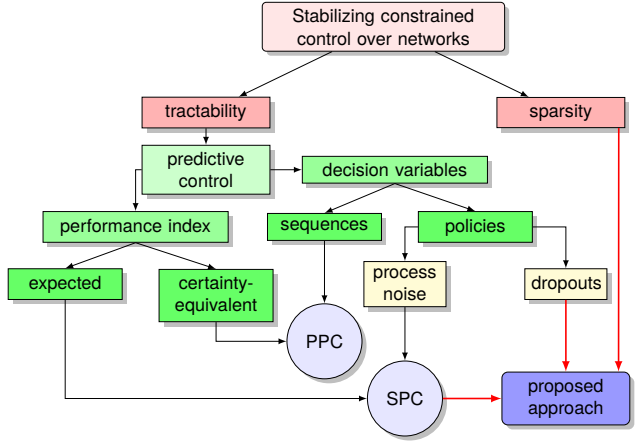

 Predictive control techniques provide tractable solutions to constrained control problems by minimizing some suitably chosen performance index over a finite temporal horizon via an iterative procedure. Based on the choice of the performance index, in context of stochastic systems, these techniques are classified as certainty-equivalent (CE) and stochastic; see Fig.\ \ref{Fig:workflow}. 
CE approaches do not take advantage of the available statistics of the uncertainties; here only the nominal plant model is considered and the control selection procedure is over open loop input sequences. CE techniques are typically implemented over networks with help of a buffer and a smart actuator; the technique is commonly known as packetised predictive control (PPC) \cite{quevedo-12}. In PPC, the time stamped sequences containing the future values of the control are transmitted at each time instant, and the successfully received sequences are saved in a buffer at the actuator. In case of dropouts, the most recent value of the control taken from the buffer is applied to the plant. 
PPC, in this way, compensates the effect of dropouts, but the controller is itself deterministic, i.e., the performance index does not incorporate the effect of unreliable communication and additive process noise. Thus, it is quite intuitive, and has been argued in \cite{ref:quevedo-15} with help of numerical experiments, that a suitably chosen stochastic performance index compensating the effect of uncertainty propagation can outperform PPC. 
	
\par Stochastic approaches incorporate the effects of uncertainties in predicted performance by considering the expected value of the cost per sample path in stochastic systems, and controlling with the help of policies as opposed to open-loop sequences. Typically, such policies are parametrized in some convenient way, and
the cost is minimized over the associated set of decision variables.\footnote{Notice that the optimization over open-loop input sequences does not give optimal performance in the presence of uncertainties \cite[pp. 13-14]{kumar1986stochastic}, and therefore, optimization over feedback policies is preferred for stochastic systems.} It is well known that feedback of past additive disturbances leads to convex problems, whereas
the state feedback approach leads to non-convexity in the set of decision variables \cite{goulart-06}. In order to obtain a convex set of feasible decision variables, disturbance feedback approaches have been studied extensively \cite{garstka1974decision,guslitser2002uncertainty,ben2004adjustable,goulart-06,ref:Lofberg-03,van2002conic}. 
To satisfy hard bounds on the control, saturated values of past disturbances are used in \cite{hokayem2009stochastic}. This saturated disturbance feedback policy is applied to networked systems with sufficient control authority in \cite{ref:amin-10} and was later generalized to any positive bound on the control in our work \cite{ref:PDQ-15}. We demonstrated in our recent conference contribution \cite{prabhatNOLCOS2016} that in the absence of additive noise, the parametrization relative to past dropouts also leads to convex problems and outperforms approaches that merely minimize over open loop input sequences. This suggests that a parametrization relative to \emph{both} past dropouts and past disturbances leads to an improved class of feedback policies. 

\par Stochastic predictive control for networked systems is based on a suitable choice of the cost function and the class of control policies, a protocol to decide what the controller will transmit and what the actuator will do. In our previous contributions \cite{ref:PDQ-15, prabhat2016, ref:quevedo-15}, we systematically developed a class of stochastic predictive control techniques for networked systems. We proposed transmission protocols in \cite{ref:PDQ-15} to answer what the controller should transmit and what the actuator should do under the class of feedback policies adopted from \cite{hokayem2009stochastic}. Stochastic approaches proposed so far \cite{ref:PDQ-15,hokayem2009stochastic,prabhat2016,ref:quevedo-15,RE-SPC} neither consider communication effects in feedback policies nor generate sparse control vectors. In this article, we propose an affine dropout and saturated disturbance feedback policy for stochastic systems controlled over unreliable and rate limited channels. Here, going beyond our earlier works, we focus on control-communication co-design by employing a sparsity promoting optimization program. We utilize the ideas of compressed sensing \cite{elad} as in \cite{bhattacharya11sparsity, nagahara2014sparse}. In \cite{bhattacharya11sparsity}, a sparsity based feedback system for the nominal plant model is presented and in \cite{nagahara2014sparse} sparse controls are designed for networked systems in absence of process noise, by exploiting a sparsity promoting regularization term, namely the $\ell_1$-norm of the control vector. 
In the present work, we employ the mixed induced $\ell_1 / \ell_{\infty}$ norm for the regularization term in presence of the feedback policy. To the best of our knowledge, this is the first work where the effects of both the process noise and the dropouts are considered in a feedback policy, sparsity in control is promoted, and stochastic stability is guaranteed.  
\par Our main contributions in this article are as follows:
\begin{itemize}[leftmargin = *]
\item We propose a policy affine in past dropouts and saturated disturbances for a finite horizon optimal control problem. The resulting problem is shown to be convex and therefore numerically tractable.
\item Stability constraints are incorporated into the underlying optimal control problem. For \emph{any} positive bound on the control and for \emph{any} non-zero successful transmission probability, these constraints ensure mean square boundedness of the system states for the largest class of linear systems with disturbances that are currently known to be stabilizable with bounded controls.
\item We introduce a regularization term in the objective function of the underlying optimal control problem to promote sparsity in time of the applied controls. Sparsity of the control commands in time is useful to reduce communication through shared channels, and increases the relaxation time for the actuator.
\item The objective function design is capable to incorporate the effects of communication channel and also control policy. This takes into account all sources of randomness in the considered networked control system, see Fig.\ \ref{Fig:blockdia}. 
\end{itemize}     
  
\par This article exposes as follows: In \secref{s:problem setup} we establish notation and definitions of the plant and its properties. In \secref{s:optimization based design} we present elementary aspects of constrained optimal control problems for stochastic systems. Our proposed class of feedback policies is presented in \secref{s:policy}. We have introduced the sparsity promoting optimal control problem in \secref{s:sparsity}. Implementation of the stabilizing feedback policy over networks is discussed in \secref{s:transmission protocol} and the computational aspects in \secref{s:computational aspects}. In \secref{s:stability}, we discuss stability issues and present stability constraints for the proposed control algorithm. We validate our results with help of numerical experiments in \secref{s:DFsimulation}. We conclude in \secref{s:epilogue}. The proofs of our main results are documented in appendix \secref{s:proofs} in consolidated manner.  

\par Our notations are standard. We let $\R$ denote the real numbers and $\Nz$ denote the positive integers. The set of the non-negative reals and non-negative integers are denoted by $\R_{\geq 0}$ and \(\Nz_0\), respectively. For any sequence \((s_n)_{n\in\Nz_0}\) taking values in some Euclidean space, we denote by \(s_{n:k}\) the vector \(\pmat{s_n\transp & s_{n+1}\transp & \cdots & s_{n+k-1}\transp}\transp\), \(k\in\Nz\). 
The notation \(\EE_z[\cdot]\) stands for the conditional expectation with given initial condition \(z\). 
For a given vector \(V\), its $i^{th}$ component is denoted by $V^{(i)}$. Similarly, $M^{(i,:)}$ and $M^{(:,i)}$ denote the $i^{th}$ row and $i^{th}$ column of a given matrix $M$, respectively. The vectors of length $k$ of all $1$'s and all $0$'s are denoted by $\ones_{k}$ and $0_k$, respectively. Similarly, \(I_d\) is the \(d\times d\) identity matrix and \(0_{r\times q} \) is the \( r\times q\) matrix with 0 entries. We use $\zeros$ in place of $0_{r \times q}$ if the dimension of the matrix is clear from the context. Inner products on Euclidean spaces are denoted by $\inprod{v}{w} \Let v \transp w$. 
\section{Plant Model and Network-System Architecture}
	\label{s:problem setup}

Throughout this work, we shall focus on a single-loop multi-variable control architecture with an unreliable input communication channel.
		Consider a linear time-invariant control system with additive process noise and controlled over an erasure channel characterised by (see Fig.\ \ref{Fig:blockdia})
		\begin{equation}
		\label{e:system}
			\st_{t+1} = \A \st_t + \B \control^a_t + \wnoise_t,\quad \st_0 = \stinit,
		\end{equation}
		where 
		\begin{enumerate}[label=(\eqref{e:system}-\alph*), leftmargin=*, widest=b]
			\item the state \(\st_t \in\R^d\), \(\stinit\in\R^d\) is a given vector. The system matrix \(A\in\R^{d\times d} \), the control matrix \( B\in\R^{d\times m}\) are given.  
\item At time $t \in \Nz_0$ the controller transmits control information $\mathcal{C}_t$, which passes through an erasure channel; see Fig.~\ref{Fig:blockdia}. This control information may not be the actual control signal. Depending upon the class of control policies and transmission protocols (discussed below in \secref{s:transmission protocol}), the control information $\mathcal{C}_t$ is used to construct the control signal by the actuator. The control $\control_t^a$ applied to the plant at time $t$ depends on the transmitted control information $\mathcal{C}_t$ and the dropouts that have occurred in the control channel till time $t$. The control $\control_t^a$ that can be delivered by the actuator takes values in the set 
				\begin{equation}\label{e:controlset}
					\controlset\Let\{v\in\R^m\mid \norm{v}_\infty \le U_{\max}\}, \text{ for each } t.
				\end{equation} 
\item The sequence \((\cnoise_t)_{t\in\Nz_0}\) is a sequence of i.i.d. Bernoulli \(\{0, 1\}\) random variables with $\EE[\cnoise_t] = p$, where $0 < p \leq 1$. 
			\item \((\wnoise_t)_{t\in\Nz_0}\) is a sequence of i.i.d. zero mean \ random vectors taking values in \(\R^d\), it is independent of \((\cnoise_t)_{t\in\Nz_0}\), and $\wnoise_t$ is symmetrically distributed around the origin for each $t$.
			\item At each \(t\) the state \(\st_t\) is measured perfectly.
			\item The communication channel between the sensors and the controller is noiseless.
            \item The acknowledgements (ACK) of the successful transmission of the control information $\mathcal{C}_t$ is causally available at the controller. 
		\end{enumerate}
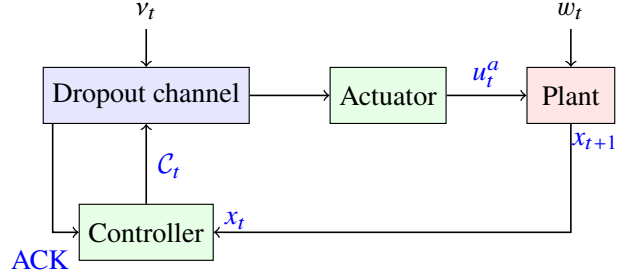
\begin{figure}[t]
\centering
\begin{adjustbox}{width = \columnwidth}
\begin{tikzpicture}
\tikzstyle{pinstyle} = [pin edge={to-,line width=0.2mm,black}]	
\tikzstyle{block} = [draw, fill=blue!10, rectangle, 
					    minimum height=2em, minimum width=1]
\tikzstyle{blockgreen} = [draw, fill=green!10, rectangle, 
					    minimum height=2em, minimum width=1cm] 
\tikzstyle{blockred} = [draw, fill=red!10, rectangle, 
					    minimum height=2em, minimum width=1.1cm] 					
\tikzstyle{sum} = [draw, circle]					     
					
\node[coordinate] (0) at (0,0) {};
\node [block, right= 0 cm of 0 ,draw, pin={[pinstyle]above:$\cnoise_t$}] (ControlChannel) {Dropout channel};
\node [blockgreen, right= 1cm of ControlChannel] (actuator) {Actuator};	
\node [blockred, right= 1cm of actuator, pin={[pinstyle]above:$\wnoise_t$}] (plant) {Plant};		
\node [blockgreen, below= 1cm of ControlChannel] (controller) {Controller};
\node[anchor=north, below= 0.7cm of plant.50, blue] {$\st_{t+1}$};
\node[anchor=north, above= 0.2cm of controller.50, blue] {$\mathcal{C}_t$};
\node[anchor=north, right= 0.2cm of actuator.20, blue] {$\control_t^a$};
\node[anchor=north, left = 0 cm of controller.202, blue] {ACK};
\node[anchor=north, right= 0 cm of controller.10, blue] {$\st_t$};

\draw[line width=0.2mm ,->] (ControlChannel.east) -- (actuator.west);
\draw[line width=0.2mm ,->] (actuator.east) -- (plant.west);
\draw[line width=0.2mm ,->] (plant.south) |- (controller.east);
\draw[line width=0.2mm ,->] (controller.north) -- (ControlChannel.south);
\draw[line width=0.2mm ,<-] (controller.west) -| (ControlChannel.197);
\end{tikzpicture}
\end{adjustbox}
\caption{Control using an unreliable channel between controller
    and actuator with causal availability of acknowledgements of successful receipt. }
\label{Fig:blockdia}
\end{figure}	
\revised{
\begin{remark}
\rm{
In many practical application the admissible control set is of the form
\begin{equation*}\label{e:ControlSet}
\controlset^{\prime} \Let \Biggl\{ v \in \R^m \Bigg\mid \abs{v^{(i)}} \leq U_i \text{ for } i =1, \ldots, m \Biggr\},
\end{equation*}  
for not necessarily equal values $U_i$.
The admissible control set $\controlset^{\prime}$ can be transformed easily into $\controlset$ as defined in \eqref{e:controlset} as follows: 
Let us define $\beta_i = \frac{U_i}{U_{\max}}$, for $i =1, \ldots, m$, and substitute the scalars $(\control_t^a)^{(i)} \Let \beta_i(v_t^a)^{(i)}$ into  \eqref{e:system} to obtain 
\begin{equation*}
	\st_{t+1} = \A \st_t + \sum_{i=1}^m \B^{(:,i)} \beta_i (v^a_t)^{(i)} + \wnoise_t.
\end{equation*}
Letting $\tilde{B}^{(:,i)} \Let \beta_i \B^{(:,i)} $, the dynamics \eqref{e:system} becomes
\begin{equation}\label{e:modifiedSystem}
\st_{t+1} = \A \st_t + \tilde{\B} v^a_t + \wnoise_t,  
\end{equation}
where $v^a_t \in \controlset$ and $\control_t^a \in \controlset^{\prime}$, by construction.
}
\end{remark}
\begin{remark}
\rm{
Recursive feasibility of stochastic predictive control techniques under state constraints is challenging whenever the additive process noise is unbounded \cite{primbs2009stochastic}. The inclusion of state constraints within our framework can be investigated in the following ways:
\begin{enumerate}[leftmargin = *]
\item Recursively feasible and stabilizing stochastic predictive control algorithm by partitioning the state space has been reported in \cite{RE-SPC}. This idea can be easily employed to incorporate state constraints in the setting of our present work. In particular, at each optimization step we check feasibility of the state constraints; if they are feasible, then we apply the controller according to our present work, otherwise apply a globally feasible recovery strategy. See \cite{RE-SPC} for details. 
\item In many applications state constraints are described by chance constraints, which are non-convex in general, and convex approximations of chance constrains have been rigorously investigated in, e.g., \cite{ben2009approximation, hong2011Convex,cinquemani-11}. An approach to approximating joint state
chance constraints by representing them as a collection of individual
chance constraints has appeared in \cite{paulson2015receding} for Schur stable systems. These individual chance constraints are further relaxed by using the Cantelli-Chebyshev inequality, and the approximated chance constrains are further softened by the exact penalty function method.
\item In some applications integrated chance constraints on the states of the form $\EE_{\st_t} \left[ \inprod{\st_t}{\mathds{S} \st_t} + \mathds{L}\transp \st_t  \right] \leq \alpha_t$, where $\mathds{S} = \mathds{S} \transp \succeq 0$ and $\alpha_t > 0$ are important; then the approach presented in \cite[Algorithm 2]{ref:Hokayem-12} can be employed. 
\item State constraints can also be studied under the framework of concentration of measure inequalities \cite{hokayem2013chance}.
\end{enumerate}
}
\end{remark}
\begin{remark}
\rm{
In this article we consider a lossy channel between the controller and the actuator but not between the sensor and the controller. Our setup caters to systems where the sensor channels have higher SNR, or guaranteed bandwidth. Examples of such systems include multi-agent systems where state-information is captured by cameras and control signals are transmitted through a wireless network, and also networks of air-borne wind energy (ABWE) systems. The extension of the present approach for a lossy channel between the sensor and the controller is non-trivial due to stochastic boundedness issues in Kalman filtering with intermittent observations \cite{intermittent, quevedo2013TAC03}. 
}
\end{remark}
\begin{remark}
\rm{
Networked systems are affected by both time delays and dropouts. To avoid book-keeping and for simplicity, here we consider delayed packets as lost packets. An explicit study of packet delays will be challenging in the setting of the present article, cf., \cite{quevedo_jurado_TAC}. 
}
\end{remark}
\begin{remark}
\rm{
We have assumed causal availability of the acknowledgements (ACK) of the successful transmissions of the control commands. For communication over one-hop links, it is a standard practice to
assume that the receiver sends an ACK signal to the transmitter for each correctly
received packet, and a negative acknowledgement (NACK) for erroneous reception. It is also reasonable to assume that the ACK and NACK signals are received error free \cite[page 207]{kuroseNetworking}.
}
\end{remark}
}
\section{Optimization Based Control Design}\label{s:optimization based design}
Let symmetric and non-negative definite matrices \(Q, Q_f\in\R^{d\times d}\) and a symmetric and positive definite matrix \(R\in\R^{m\times m}\) be given. We define a standard quadratic cost-per-stage function \(\costps:\R^d\times\controlset\lra \R_{\geq 0} \) and a final cost function \(\costfinal:\R^d\lra \R_{\geq 0} \) by \[ \costps(z, v) \Let \inprod{z}{Qz} + \inprod{v}{Rv} \text{ and } \costfinal(z) = \inprod{z}{Q_f z},\] respectively. Fix an optimization horizon \(N\in\Nz\) and consider the objective function at time \(t\) given the state \(\st_t\):
		\begin{equation}
		\label{e:objfn}
			V_t \Let \EE_{\st_t}\biggl[\sum_{k=0}^{N-1} \costps(\st_{t+k}, \control^a_{t+k}) + \costfinal(\st_{t+N}) \biggr].
		\end{equation}
		In this setting the cost function \(V_t\), intuitively, considers the control effort that occurs at the actuator end, not just the computed control. This makes more sense than considering the effort with respect to the control commands since the commands may be dropped by the erasure channel.\footnote{ In \secref{s:sparsity} we show how to design the cost function to obtain sparse control transmissions.} At time instant \(t\), we are interested in minimizing the objective function \(V_t\) over a class of causal history-dependent feedback strategies \(\Pi\) formally defined by 
		\[
\control_{t+\ell} =
\begin{cases}		
 \pi_{t+\ell}(\st_t, \cdots, \st_{t+\ell}, \cnoise_t, \cdots, \cnoise_{t+\ell-1}) \quad \text{ for } \ell = 1, \cdots, N-1 \\
 \pi_{t}(\st_t) \quad \text{ for } \ell = 0
\end{cases}			
		\]
		while satisfying \(\control_{t+\ell} \in\controlset\) for $\ell = 0,1, \ldots, N-1$ for each $t$. Recall that a control strategy or policy \(\pi\) is a sequence \((\pi_0, \ldots, \pi_t, \ldots)\) of Borel measurable maps \(\pi_t:\R^d\lra\controlset\) \cite[\S 2.1]{surveyMDP}. Policies of finite length \((\pi_t, \pi_{t+1}, \ldots, \pi_{t+N-1})\) for some \(N\in\Nz\) will be denoted by \(\pi_{t:N}\) in the sequel. We refer the readers to \cite[\S II]{dualEffect} for a succinct discussion on the classes of policies. We revisit the control policy classes in \secref{s:policy} again and refer to them as \emph{history dependent feedback policies} in connection with the above references.
		\par
\begin{figure}[t]
\begin{adjustbox}{width = \columnwidth}
        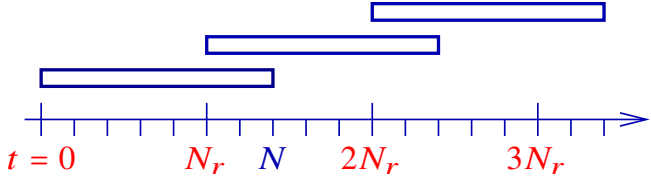
\end{adjustbox}
\caption{Receding horizon control strategy: At $t=0$, a control sequence of length $N$ is computed and the first $N_r$ controls of that sequence are applied to the plant. This process repeats after every $N_r$ time steps.}
\label{Fig:recalculation}
\end{figure}

 The \emph{receding horizon control strategy} for a given \emph{recalculation interval} \(N_r \in\{1, \ldots, N\}\) consists of successive applications of the following steps (see Fig.\ \ref{Fig:recalculation}):
		\begin{enumerate}[label=(\roman*), leftmargin=*, widest=iii]
			\item measure the state \(\st_t\) and determine an admissible feedback policy \(\pi^\star_{t:N}\in\Pi\), \label{step:compute}
			\item apply the first \(N_r \) elements \(\pi^\star_{t:N_r-1}\) of this policy,
			\item increase \(t\) to \(t+ N_r \) and return to step (i).
		\end{enumerate}
In order to avoid step \ref{step:compute} at each time instant one should use $N_r > 1$ by following the approach presented in this article. However, the presented approach is valid for $N_r=1$ as well. In \secref{s:stability} Remark \ref{rem:recalculation} we present more clarification on the choice of $N_r$. 

		The states, controls, and noise over one optimization horizon of length \(N\) admit the following description under an \textsl{unreliable control channel}: 
		\begin{equation}
		\label{e:augmented}
			\st_{t:N+1} = \calA \st_t + \calB \control^a_{t:N} + \calD \wnoise_{t:N},
		\end{equation}
		where $\calA$, $\calB$ and $\calD$ are standard matrices of appropriate dimensions.
We define two block diagonal matrices $\calQ$ and $\calR$ of appropriate dimensions
to state the following preliminary versions of the optimal control problem underlying the receding horizon control technique:
		\begin{equation}\label{e:opt problem}
		\begin{aligned}
			& \minimize_{\pi_{t:N}}	&& \EE_{\st_t}\bigl[ \inprod{\st_{t:N+1}}{\calQ \st_{t:N+1}} + \inprod{\control^a_{t:N}}{\calR\control^a_{t:N}} \bigr]\\
			& \sbjto	&& 
				\begin{cases}
					\text{constraint } \eqref{e:augmented},\\
					\control_t^a \in\controlset\;\text{for all }t,\\
					\pi_{t:N}  \text{ in a class of policies}.
				\end{cases}
		\end{aligned}
		\end{equation}
To account for the presence of control channel noise, the optimal control problem \eqref{e:opt problem} will be refined by an appropriate selection of policies to be addressed in \secref{s:policy}. In particular, we shall employ feedback from the process noise and the control channel dropouts in our policies. This will lead to a modification of the optimal control problem \eqref{e:opt problem}. 
The scheme of the transmission of the parameters affects the optimization problem. Therefore, the transmission protocol must also be considered at the control synthesis stage. We shall include the effect of transmission protocol on the optimization problem with help of the stacked control vector $\control_{t:N}^a$ in \secref{s:transmission protocol}.
	\section{Control Policy Class}\label{s:policy}

We recall that, by assumption, the states are completely and exactly measured, and acknowledgment of whether a successful control transmission has occurred or not is assumed to be causally available to the controller, see Fig.\ \ref{Fig:blockdia}. Therefore, the dynamics of the actuator is known to the controller and it is possible to reconstruct the noise sequence from the sequence of observed states and control inputs with the aid of the formula
		\begin{equation}\label{e:estimation}
			\wnoise_t = \st_{t+1} - \A \st_t - \B \control^a_t,\quad t\in\Nz.
	    \end{equation}
This calculation is performed by the controller at every time \( t \in \Nz_0\).	   
In our earlier work \cite{ref:PDQ-15} we followed the approach developed in \cite{ref:ChaHokLyg-11} and employed \(N\)-history-dependent policies, which are affine in saturated values of past disturbances and an offset parameter. 
In the present work, we build on the fact that there are two sources of uncertainty, namely, dropouts and additive process noise. Accordingly, we refine the \(N\)-history-dependent policies in \cite{ref:PDQ-15} to include the control channel dropouts as: 
\begin{equation}\label{e:controlt}
			\control_{t+\ell} = \eta_{t+\ell} + \sum_{i=0}^{\ell-1} \Bigl( \theta_{\ell,t+i} \ee_{i+1}(\wnoise_{t+i}) + \lambda_{\ell,t+i} \cnoise_{t+i} \Bigr),
		\end{equation}		
for \(\ell =0,1, \ldots, N-1\). In \eqref{e:controlt}, $\theta_{\ell, t+i}$'s and $\lambda_{\ell,t+i}$'s are the policy parameters.		
		 The control vector \(\control_{t:N}\) admits the compact representation:
		\begin{equation}\label{e:policy}
			\control_{t:N}	 
		 \teL \offset_t + \gain_t \ee(\wnoise_{t:N-1}) + \weight_t \cnoise_{t:N-1}. 
		\end{equation}
		In \eqref{e:policy}, the component-wise saturation function $\ee \Let \pmat{\ee_1 \transp &  \ee_{2} \transp & \cdots & \ee_{N-1}\transp}\transp$, the offset vector \(\offset_t\in\R^{m N}\) and the gain matrices \(\gain_t , \weight_t \) are  strictly lower block triangular matrices
		\begin{equation} \label{e:gain}
		\begin{aligned}
			\gain_t &= \pmat{0 & 0 & \cdots & 0 & 0\\ \theta_{1, t} & 0 & \cdots & 0 & 0\\ \theta_{2, t} & \theta_{2, t+1}  & \cdots & 0 & 0\\ \vdots & \vdots & \vdots & \vdots & \vdots\\ \theta_{N-1, t} & \theta_{N-1, t+1} & \cdots & \theta_{N-1, t+N-3} & \theta_{N-1, t+N-2}},\\
			\weight_t &= \pmat{0 & 0 & \cdots & 0 & 0\\ \lambda_{1, t} & 0 & \cdots & 0 & 0\\ \lambda_{2, t} & \lambda_{2, t+1}  & \cdots & 0 & 0\\ \vdots & \vdots & \vdots & \vdots & \vdots\\ \lambda_{N-1, t} & \lambda_{N-1, t+1} & \cdots & \lambda_{N-1, t+N-3} & \lambda_{N-1, t+N-2}},
			\end{aligned}
		\end{equation}
		with each \(\theta_{k, \ell} \in \R^{m\times d}, \lambda_{k, \ell} \in \R^{m}\), and $\norm{\ee(\wnoise_{t:N-1})}_{\infty} \leq \varphi_{\max}$. 
We assume that the component-wise saturation function \(\ee\) is symmetric about the origin. For instance, we can employ either standard saturation, piecewise linear, or sigmoidal functions \cite[\S III]{ref:ChaHokLyg-11}. The saturation functions $\ee_i$'s can be chosen different for each $i$ when the variance changes within a prediction horizon. Also, different bounds for each component of $\ee_i$ can be chosen when the components of the disturbance at a given time have different variance, see also \cite{korda2012}.  
\section{Sparsity Promoting Optimization-based Control}\label{s:sparsity}
Since in networked systems rate limited channels are shared among various devices, sparse controls are also desirable and tractability is essential for the implementation. 
The solution of the optimal control problem \eqref{e:opt problem} at time $t$ under the class of control policies \eqref{e:policy} gives offset vector $\offset_t$, and gain matrices $\gain_t$ and $\weight_t$. A control vector $\control_{t+\ell}$ for $\ell \in \{ 0,\cdots,N-1 \}$ is equal to zero when the $\ell^{th}$ block of the matrix $F_t \Let \pmat{\offset_t & \gain_t & \weight_t}$ has all entries equal to zero: \[(F_t)_{\ell m: (\ell+1)m-1} = 0_{m \times (d+1)(N-1)+1}. \]
Generally, a sparsity promoting convex optimal control problem employs $\ell_1$ norm as a regularization term when decision variable is a vector \cite{nagahara2014sparse}. When decision variable is in form of a matrix, $\ell_1 / \ell_{\infty}$ mixed induced norm is used as a regularization term \cite{sparsityNorm}. In order to get sparse matrix $F_t$ with zero matrices of the dimension $m \times (d+1)(N-1)+1$, we use $\ell_1 / \ell_{\infty}$ mixed induced norm of $F_t$ when $m=1$. For $m>1$, we construct another matrix $\hat{F}_t$ by taking transpose of each $m\times 1$ dimensional block of $F_t$. By doing so we obtain a matrix $\hat{F_t}$ of dimension $N \times (m(d+1)+1)$. We define the regularization term 
\begin{equation} \label{e:regularization term}
\regular_t(\offset_t,\gain_t,\weight_t) \Let \sum_{i=1}^{N} \norm{\hat{F}_t^{(i,:)}}_{\infty} 
\end{equation}
For $\mu \geq 0$ we have the following optimal control problem
	\begin{equation}\label{e:opt problem sparse}
		\begin{aligned}
			& \minimize_{\pi_{t:N}}	&& \EE_{\st_t}\bigl[ \inprod{\st_{t:N+1}}{\calQ \st_{t:N+1}} + \inprod{\control^a_{t:N}}{\calR\control^a_{t:N}} \bigr]\\
			& ~ && + \mu \regular_t(\offset_t,\gain_t,\weight_t) \\ 
			& \sbjto	&& 
				\begin{cases}
					\text{constraint } \eqref{e:augmented},\\
					\control_t^a \in\controlset\;\text{for all }t,\\
					\pi_{t:N}  \text{ in a class of policies}.
				\end{cases}
		\end{aligned}
		\end{equation}
\section{Implementation of Control Policy over Networks}\label{s:transmission protocol}
The class of control policies proposed in \secref{s:policy} is implemented over networks with the help of transmission protocols. The parameters of the control policy can be transmitted in several different ways through the control channel.\footnote{Please see \cite{ref:PDQ-15} for a detailed discussion on transmission protocols.} For demonstration purposes, in this article we consider the following two transmission protocols:

\subsection{Sequential Transmission of Control Values} 

		\begin{enumerate}[label={\rm (TP\arabic*)}, leftmargin=*, widest=3, align=left, start=1]
			\item \label{a:seq} Solely the control values $\control_{t+\ell}$ are computed by the controller and transmitted to the actuator at each instant $t+\ell$, ~ $\ell \in \{0, \ldots, ~ N_r -1 \}, ~ t = K N_r, ~ K \in \Nz_0$.
		\end{enumerate}

The control value transmitted at time $t+\ell$ is affected by dropout $\cnoise_{t+\ell}$. Hence, $\control^a_{t:N}$ is given by   
		\begin{equation}
		\label{e:policyseq}
			\control^a_{t:N} \Let \calS \Bigl(\offset_t + \gain_t\ee(\wnoise_{t:N-1}) + \weight_t\cnoise_{t:N-1} \Bigr), 
		\end{equation}	
			
		where $\calS \Let \blkdiag\bigl(I_{m}\cnoise_{t}, \cdots, I_{m} \cnoise_{t + N_r -1}, I_{m(N-N_r)} \bigr)$.
In the above transmission protocol, if the control packet is lost at some time instant, zero control is applied to the plant at that instant.
Note that \ref{a:seq} does not require any storage or computational facility at the actuator end.

\subsection{Sequential transmission of control along with repetitive transmission of remaining offset components}
\begin{enumerate}[label={\rm (TP\arabic*)}, leftmargin=*, widest=3, align=left, start=2]
\item \label{a:repetitive} The control value $\control_{t+\ell}$  is transmitted at each instant $t+\ell$, $\ell \in \{0, \ldots, N_r -1 \}, ~ t = K N_r, ~ K \in \Nz_0$. In addition, the remaining blocks of the current offset vector $\offset_t$ are transmitted at each step until the first successful transmission. 
\end{enumerate}
		 
According to \ref{a:repetitive}, control values are transmitted at each instant as in \ref{a:seq}. To mitigate the effects of dropouts, those components of the burst $(\offset_t)_{1:mN_r}$ that may become useful at future instants are also transmitted repetitively until one packet is successfully received at the actuator. Successfully received packets are stored in the buffer, so that in the case of packet dropout, the corresponding offset block is used. 
For $\ell = 0,1, \cdots, N_r -2$, the control values along with the burst of the remaining offset components transmitted at time $t+\ell$ are affected by the dropout at the same time instant. The plant noise $\wnoise_{t+ \ell}$ is recovered at the controller by \eqref{e:estimation} correctly with the help of causally available acknowledgements. The plant input sequence using \ref{a:repetitive} can therefore be represented in compact form as:   
        \begin{equation}\label{e:policyrepetitive}
			\control^a_{t:N} \Let \mathcal{G}\offset_t + \calS\left(\gain_t\ee(\wnoise_{t:N-1}) + \weight_t \cnoise_{t:N-1} \right), \\
		\end{equation}
		where the matrix $\calG$ has $(N \times (N-1))$ blocks in total, each of dimension $m \times m$. For $i = 1,\cdots,N$ and $j = 1,\cdots,N-1$, the matrix $\calG$ can be given in terms of the blocks $\calG_b^{(i,j)}$ each of dimension $m \times m$ as follows: 		 
			\begin{align}
				\calG_b^{(i,j)} \Let
		\begin{cases}
		     \rho_{t+i-1} I_m & \text{ if } i=j \leq N_r ,\\ 
			I_{m} & \text{ if } i =j > N_r ,\\
			0_m \quad  & \text{ otherwise,}			
		\end{cases}		
				\end{align}
		where $\rho_t = \cnoise_t$, $\rho_{t+\ell} = \rho_{t+\ell-1} + \left( \prod_{s = 0}^{\ell -1} (1-\cnoise_{t+s}) \right) \cnoise_{t+\ell}  $, and $\gain_t, \weight_t$ and $\calS$ are given in \eqref{e:gain} and \eqref{e:policyseq},  respectively; and $\ee(\wnoise_{t:N-1})$ is as defined in \eqref{e:policy}. The term $\rho_{t+\ell}$ captures the effect of \ref{a:repetitive}. 
		Note that \ref{a:repetitive} requires storage facility at the actuator, but no advanced computation capacity.

\section{Computational Aspects} \label{s:computational aspects}
The optimal control problem \eqref{e:opt problem sparse} must be solved periodically online. In order to write the objective function of \eqref{e:opt problem sparse} in terms of variance and covariance matrices, which can be computed off-line to reduce the on-line computational burden, we manipulate the decision variable $\weight_t$ as discussed below. Let us consider \[ \mathfrak{C}_{\ell} \Let \pmat{0_{m\ell \times m(N-\ell)} \\ I_{m(N-\ell) \times m(N-\ell)} } \]  and let $\bar{\weight}_t^{(:,\ell)} \Let \mathfrak{C}_{\ell} \transp \weight_t^{(:,\ell)} $ contains non-zero entries of the $\ell^{th}$ column vector of $\weight_t$ which is obtained by removing the first $\ell m$ entries of $\ell^{th}$ column of $\weight_t$. Let us define 
\begin{equation}\label{e:manipulated}
\Xi_t \Let \pmat{\offset_t \\ \tilde{\weight}_t }, 
\end{equation}
where $\tilde{\weight}_t \Let \pmat{(\bar{\weight}_t^{(:,1)})\transp & (\bar{\weight}_t^{(:,2)})\transp & \cdots & (\bar{\weight}_t^{(:,N-1)})\transp }\transp $ consists of concatenated column vectors of $\weight_t$ with non-zero entries. We have the following theorems: 
\begin{thm}
		\label{th:DF_rep}
			Consider the control system \eqref{e:system} 
			under the transmission protocol \ref{a:repetitive}. Then,
			 for every $t=0, N_r, 2N_r,\ldots$, the optimization problem \eqref{e:opt problem sparse} over the class of control policies \eqref{e:controlt} is convex, feasible, and can be rewritten as the following convex quadratic program:
\begin{equation}\label{e:programsingle}
\begin{aligned}
\minimize_{\offset_t ,\gain_t, \weight_t }\;	 &\;  \inprod{\Xi_t}{ \calL \Xi_t } + \inprod{\calM \calA \st_t}{  \Xi_t } +  \inprod{\calA \st_t}{\calQ \calA \st_t} \\
& \quad+ 2\trace(\gain_t\transp\mu_{\calS}\transp\calB\transp\calQ\calD\Sigma_{\ee}^{\prime}) + \trace(\calD \transp \calQ \calD \Sigma_{\wnoise}) \\
&\quad + \trace(\gain_t \transp \Sigma_{\calS} \gain_t \Sigma_{\ee}) + \mu\regular_t(\offset_t,\gain_t,\weight_t)
\end{aligned}
\end{equation}
 \quad \( \sbjto \) :
\begin{align}\label{e:decisionboundsingle}
\abs{\offset_t^{(i)} + \dfrac{1}{2} \weight_t^{(i,:)}\ones_{(N-1)}} + \dfrac{1}{2}\norm{\weight_t^{(i,:)}}_1 + \norm{\gain_t^{(i,:)}}_1\varphi_{\max} & \leq U_{\max},
\end{align}
for  $ i = 1, \cdots, mN$, \\
%
where 
\begin{equation}
\label{e:calL}
\calL = \pmat{\Sigma_{\calG} & \tilde{\Sigma}_{\calS\calG_{\ell}} \\ \tilde{\Sigma}_{\calS\calG_{\ell}}\transp & \tilde{\Sigma}_{\calS_{ n \ell}} },
\end{equation}
 $\Sigma_{\calG} = \EE[\calG\transp\alpha\calG]$, $\alpha = \calB\transp\calQ\calB+\calR$, \\
 $\tilde{\Sigma}_{\calS\calG_{\ell}} = \pmat{ \bar{\Sigma}_{\calS\calG_1} & \bar{\Sigma}_{\calS\calG_2} & \cdots & \bar{\Sigma}_{\calS\calG_{(N-1)}} }$, $\bar{\Sigma}_{\calS\calG_{\ell}} = \Sigma_{\calS\calG_{\ell}} \mathfrak{C}_{\ell}$, $\Sigma_{\calS\calG_{\ell}} = \EE\Bigl[\calG \transp\alpha \calS_{\ell} \Bigr]$,  $\calS_{\ell} \Let \cnoise_{t+\ell-1} \calS$, \[ \tilde{\Sigma}_{\calS_{n \ell}} = \pmat{\bar{\Sigma}_{\calS_{11}} & \bar{\Sigma}_{\calS_{12}} & \cdots & \bar{\Sigma}_{\calS_{ 1(N-1)}}\\ \vdots & \vdots & \cdots & \vdots \\ \bar{\Sigma}_{\calS_{ (N-1) 1}} & \bar{\Sigma}_{\calS_{ (N-1)2}} & \cdots & \bar{\Sigma}_{\calS_{ (N-1)(N-1)}}}, \] 
$\bar{\Sigma}_{\calS_{n\ell}} = \mathfrak{C}_n \transp \Sigma_{\calS_{n \ell}} \mathfrak{C}_{\ell}$, $\Sigma_{\calS_{ n\ell}} = \EE\Bigl[\calS_{n} \transp\alpha  \calS_{\ell} \Bigr] $,
\begin{equation}
\label{e:calM}
\calM = 2\pmat{\calQ \calB \mu_{\calG} & \calQ \calB \tilde{\mu}_{\calS_{\ell}} }\transp,
\end{equation} 
$\mu_{\calG}=\EE[\calG]$, $\tilde{\mu}_{\calS_{\ell}}= \pmat{\bar{\mu}_{\calS_1} & \bar{\mu}_{\calS_2} & \cdots & \bar{\mu}_{\calS_{(N-1)}} }$,  
$\bar{\mu}_{\calS_{\ell}} = \mu_{\calS_{\ell}} \mathfrak{C}_{\ell}$, $\mu_{\calS_{\ell}} = \EE\left[\calS_{\ell}\right]$,     $\Sigma_e^{\prime} \Let \EE[\wnoise_{t:N} \ee(\wnoise_{t:N-1})\transp]$,  $\Sigma_W \Let \EE[\wnoise_{t:N}\wnoise_{t:N}\transp]$,        
$\Sigma_e \Let \EE[\ee(\wnoise_{t:N-1})\ee(\wnoise_{t:N-1})\transp]$, $\regular_t(\offset_t, \gain_t, \weight_t)$ is defined in \eqref{e:regularization term} and $\Xi_t$ in \eqref{e:manipulated}.
\end{thm}				
\begin{thm}
\label{th:DF_seq}
			Consider the control system \eqref{e:system} under the transmission protocol \ref{a:seq}. Then,
			 for every $t=0, N_r, 2N_r,\ldots$, the optimization problem \eqref{e:opt problem sparse} over the class of control policies \eqref{e:controlt} is convex, feasible, and can be rewritten as the convex quadratic program as in Theorem \ref{th:DF_rep} by substituting $\calS$ in place of $\calG$ in $\calL$ and $\calM$ in \eqref{e:calL} and \eqref{e:calM}, respectively.
\end{thm}
\begin{remark}
\rm{The covariance matrices $\calL$, $\calM$, $\Sigma_{\ee}$, $\Sigma_W$ and $\Sigma_{\ee^{\prime}}$ that are required to solve the optimization problem were computed empirically via classical Monte Carlo methods \cite{ref:robert-13} using $10^6$ i.i.d. samples. Computations for determining our policy are carried out in the MATLAB-based software package YALMIP \cite{ref:lofberg-04} and are solved using SDPT3-4.0 \cite{ref:toh-06}. 
}
\end{remark}
\revised{
\begin{remark}
\rm{
Since the matrix pair $(\A,\B)$ is not affected by uncertainties, here we present a tractable solution by separating all sources of uncertainties in terms of the variance and the covariance matrices $\calL$, $\calM$, $\Sigma_{\ee}$, $\Sigma_W$. When there is an implicit mapping between uncertainties and the states, the polynomial chaos framework is often employed to obtain an approximate explicit mapping between the uncertainties and the states \cite{paulson2015stability}; this is especially useful for nonlinear systems than for linear ones with noisy system and controller models. 
}
\end{remark}
\begin{remark}
\rm{
Our approach presents a quadratic program with number of decision variables $\tilde{N} = mN(1+ \frac{N-1}{2}(d+1))$, which gives an overall complexity of our program as $O(\tilde{N}^3)$ \cite{cormen2009introduction}.
Moreover, as carried out in \cite[Example 6]{ref:ChaHokLyg-11} by setting elements of $\Lambda_t$ and $\Theta_t$ to $0$ except subdiagonal elements, the number of decision variables reduces to $(N-1)m(d+1)+m$.        
}
\end{remark}
}
\section{Stability Issues} \label{s:stability}	
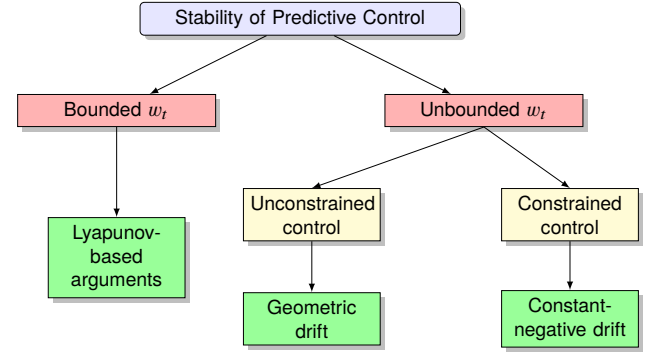
\begin{figure}
\begin{adjustbox}{width = \columnwidth}


\tikzset{
  basic/.style  = {draw, drop shadow, font=\sffamily, rectangle},
  root/.style   = {basic, rounded corners=2pt, thin, align=center,
                   fill=blue!10, text width=5cm},
  level 2/.style = {basic, thin,align=center, fill=red!30,
                   text width=3cm},
  level 3/.style = {basic, thin, align=center, fill=yellow!20, text width=2cm},
  level 4/.style = {basic, thin, align=center, fill=green!40, text width=2cm},
  }

\begin{tikzpicture}[
  level 1/.style={sibling distance=60mm},
  edge from parent/.style={->,draw},
  >=latex]

\node[root] {Stability of Predictive Control}
child {node[level 2] (bounded) {Bounded $\wnoise_t$}}
child {node[level 2] (unbounded) {Unbounded $\wnoise_t$}};

\begin{scope}[every node/.style={level 3}]
\node [below of = unbounded, xshift = -80pt, yshift = -20pt] (unboundedC) {Unconstrained control};
\node [below of = unbounded, xshift = +40pt, yshift = -20pt] (boundedC) {Constrained control};
\end{scope}

\begin{scope}[every node/.style={level 4}]
\node [below of = unboundedC, yshift = -20pt] (geometric) {Geometric drift};
\node [below of = boundedC, yshift = -20pt] (negative) {Constant-negative drift};
\node [below of = bounded, yshift = -40pt] (lyapunov) {Lyapunov-based arguments};
\end{scope}

\draw[->] (bounded.south) -- (lyapunov.north);
\draw[->] (unbounded.south) -- (boundedC.north);
\draw[->] (unbounded.south) -- (unboundedC.north);
\draw[->] (unboundedC.south) -- (geometric.north);
\draw[->] (boundedC.south) -- (negative.north);

\end{tikzpicture}
\end{adjustbox} 
\caption{In the presence of constraints on controls and the unbounded process noise, the constant-negative drift conditions are used.}
\label{fig:StabilityTree}
\end{figure}

\revised{
Stability in optimization based control techniques is typically achieved by selecting either approximate cost functions satisfying some Lyapunov based conditions, or by enforcing stability constraints in the underlying optimal control problem \cite[\S 3.8.3]{ref:mayne-00}. Both approaches are conservative in general, and ensuring good closed-loop behaviour in the presence of bounded control authority is difficult, see \cite{meyn2012markov, ref:ChaRamHokLyg-12}. We refer readers to \cite[\S IV]{ref:PDQ-15} for a detailed discussion, a summary of which is given in Fig.\ \ref{fig:StabilityTree}. 
}
\revised{
It has been demonstrated in \cite{Sussmann97} and \cite[Theorem 1.7]{ref:ChaRamHokLyg-12} that one \emph{cannot} globally stabilize an LTI system \eqref{e:system}
 by bounded control actions (even in the presence of a perfect channel) if the spectral radius of the system matrix $A$ is greater than unity. In view of this fundamental restriction, we make the following assumption:
 \begin{assumption}\label{as:lyapunov}
The system matrix $\A$ is Lyapunov stable and the matrix pair $(A,B)$ is stabilizable.\footnote{Recall that a Lyapunov stable matrix (neutral system) has its eigenvalues bounded by unit circle. Those eigenvalues located on the unit circle have equal algebraic and geometric multiplicities \cite[page 211]{stabilitybook}.}
 \end{assumption}
In the light of Assumption \ref{as:lyapunov}, the system dynamics \eqref{e:system} takes the following form: 
\begin{equation}
		\label{e:orthogonal decomposition}
		\pmat{\stortho_{t+1} \\ \st^s_{t+1}} = \pmat{\Aortho & 0\\ 0 & \Aschur}\pmat{\stortho_{t} \\ \st^s_{t}}	+ \pmat{\Bortho\\ \Bschur}\control_t + \pmat{\wnoise_t^o \\ \wnoise_t^s}, 
		\end{equation}
where \(\Aortho\in\R^{d_o\times d_o}\) is orthogonal and \(\Aschur\in\R^{d_s\times d_s}\) is Schur stable, with state dimension \(d = d_o + d_s\). Since the pair \((\A, \B)\) is assumed stabilizable, there exists a positive integer \(\reachindex\) such that the pair \((\Aortho, \Bortho)\) is reachable in \(\reachindex\)-steps. If $	\reachab_{\ell} \Let \pmat{\Aortho^{\ell-1}\Bortho & \Aortho^{\ell-2}\Bortho & \cdots & \Bortho}$, then $\rank(\reachab_{\reachindex}) = d_o$.  
This positive integer $\reachindex$ is called the \emph{reachability index} of the pair $( \Aortho, \Bortho )$ in the sequel. \\
We shall focus on mean-square boundedness of the closed-loop system \eqref{e:system}, and we recall the definition of mean square boundedness \cite[\S III.A]{chatterjee-15} and recast Lemma  1 and Lemma 2 from \cite{ref:PDQ-15} for our subsequent analysis: 
\begin{lem}{\cite[Lemma 1, Lemma 2]{ref:PDQ-15}}
\label{t:msbsingle}
Consider the orthogonal part of the system \eqref{e:system} given
as per the decomposition \eqref{e:orthogonal decomposition}. If there exists a \(\reachindex\)-history dependent policy and $\control_t \in \controlset$ is chosen such that for any given $\epsilon, r > 0$, $0 < \zeta < \frac{U_{\max}}{\sqrt{d_o}\sigma_{1}\left(\reachab_{\reachindex}^{+}\right)}$, and for any $t=0, \reachindex, 2\reachindex, ...$, $ j=1,2,\cdots ,d_o$, the following drift conditions are satisfied: 
\begin{align}
\EE_{\stortho_t} \left[ \left( (\Aortho^{\reachindex})\transp \reachab_{\reachindex}(\Aortho, \Bortho)\control_{t:\reachindex} \right)^{(j)} \right] & \leq -\zeta \nonumber \\
\text{ whenever } \left( \stortho_{t} \right)^{(j)} & \geq r + \epsilon \label{e:drift1} \\
\EE_{\stortho_t} \left[ \left( (\Aortho^{\reachindex})\transp \reachab_{\reachindex}(\Aortho, \Bortho)\control_{t:\reachindex}\right)^{(j)} \right] & \geq \zeta \nonumber \\
\text{ whenever }\left( \stortho_{t}  \right)^{(j)} & \leq -r - \epsilon, \label{e:drift2}
\end{align} 
then this policy renders the orthogonal part of the closed-loop system \eqref{e:orthogonal decomposition} mean-square bounded under both the transmission protocols \ref{a:seq} and \ref{a:repetitive} introduced in Section 6.
\end{lem} 
To achieve the drift \eqref{e:drift1} and \eqref{e:drift2}, we write the first $\reachindex$ blocks in \eqref{e:controlt} as 
\begin{equation}\label{e:decision}
\control_{t: \reachindex} \Let (\offset_t)_{1: \reachindex m} +(\weight_t)_{1: \reachindex m}\cnoise_{t:N-1} + (\gain_t)_{1:\reachindex m}\wnoise_{t:N-1},
\end{equation}
and substitute \eqref{e:decision} into \eqref{e:drift1} and \eqref{e:drift2}. Now,
for given $\epsilon, r >0$ and for every $j=1,\cdots,d_o$, we arrive at the following stability constraints, in lieu of \eqref{e:drift1} and \eqref{e:drift2}:
\begin{align}
\Bigl( (\Aortho^{\reachindex})\transp \reachab_{\reachindex}\left[(\offset_t)_{1:\reachindex m} +p(\weight_t)_{1: \reachindex m} \right] \Bigr)^{(j)} &\leq -\zeta \nonumber  \\
\text{ whenever} \left( \stortho_{t} \right)^{(j)} &\geq r + \epsilon, \label{e:decisonConstraint1single} \\
\Bigl( (\Aortho^{\reachindex})\transp \reachab_{\reachindex}\left[(\offset_t)_{1:\reachindex m} +p(\weight_t)_{1: \reachindex m} \right] \Bigr)^{(j)}  & \geq \zeta \nonumber \\
\text{ whenever} \left( \stortho_{t} \right)^{(j)} & \leq -r - \epsilon. \label{e:decisonConstraint2single}
				\end{align}
}				
Having established the above, for mean square boundedness of \eqref{e:system} in closed-loop under \ref{a:seq} and \ref{a:repetitive} we have the following results:
\begin{thm}
		\label{th:DF_rep_stbl}
			Consider the control system \eqref{e:system} under the transmission protocol \ref{a:repetitive}, and let Assumption \ref{as:lyapunov} hold.
			\begin{enumerate}[label={\rm (\roman*)}, leftmargin=*, widest=ii, align=right]
				\item For every $t=0, \reachindex, 2\reachindex,\ldots$, the optimization problem \eqref{e:opt problem sparse} over the class of control policies \eqref{e:controlt} along with the stability constraints \eqref{e:decisonConstraint1single}-\eqref{e:decisonConstraint2single} is convex, feasible, and can be rewritten as the following convex quadratic program:
\begin{equation*}
\begin{aligned}
			\minimize_{\offset_t ,\gain_t, \weight_t } \quad & \text{Objective function}  \; \eqref{e:programsingle} \\
			\sbjto \quad & \text{Constraints } \eqref{e:decisionboundsingle}, \eqref{e:decisonConstraint1single} \text{ and } \eqref{e:decisonConstraint2single}.\\
\end{aligned}
		\end{equation*}
				\item For any initial condition \(\stinit \in \R^{d}\) successive application of the control law given by the optimization problem in (i) for \( \reachindex \) steps renders the closed-loop system mean square bounded.
				\end{enumerate}
		\end{thm}
		
\begin{thm}\label{th:DF_seq_stbl}
For the control system \eqref{e:system} under the transmission protocol \ref{a:seq}, the assertions of Theorem \ref{th:DF_rep_stbl} hold by substituting $\calS$ in place of $\calG$ in $\calL$ and $\calM$ in \eqref{e:calL} and \eqref{e:calM}, respectively.
\end{thm}	
\begin{remark}
\rm{ Theorem \ref{th:DF_rep_stbl} generalizes the main result of \cite{ref:PDQ-15}. Indeed, if we set $\weight_t = 0$ in \eqref{e:policy} and and $\mu = 0$ in \eqref{e:programsingle}, we recover the main result of \cite{ref:PDQ-15}. Theorem \ref{th:DF_seq_stbl} generalizes the main result of \cite{prabhatNOLCOS2016} for the case $\Sigma_{\wnoise} = 0$. Indeed, if we set $\gain_t = 0$ in \eqref{e:policy} and $\mu = 0$ in \eqref{e:programsingle}, we recover the main result of \cite{prabhatNOLCOS2016}. 
}
\end{remark}
\begin{remark}\label{rem:recalculation}
{\rm
We need the recalculation interval $N_r$ equal to the reachability index $\reachindex$ of the orthogonal sub-system to satisfy the drift conditions \eqref{e:drift1} and \eqref{e:drift2} above.
}
\end{remark}
\section{Numerical Experiments} \label{s:DFsimulation}

  
In this section we present simulations to illustrate our results. Consider the three dimensional linear stochastic system 
\begin{equation}\label{e:example}
\st_{t+1} = \pmat{0&-0.80&-0.60\\ 0.80&-0.36&0.48\\ 0.60&0.48&-0.64}\st_t + \pmat{0.16\\0.14\\1}\control_t + \wnoise_t, 
\end{equation}
where $\wnoise_t$ is i.i.d. Gaussian of mean zero, the initial condition is $\stinit=\pmat{10&10&-10} \transp$, and the control is bounded as per $|u_t|\leq 15$. 
\par We solved a constrained finite-horizon optimal control problem corresponding to the states and the control weights \[ Q = I_3, Q_f = \pmat{12& 1& 4\\ 1&19& 2\\ 4& 2&2}, R =2 .\] We selected an optimization horizon, \(N=4\), recalculation interval \(N_r = \reachindex = 3 \) and simulated the system responses. We selected the nonlinear bounded term $\ee(\wnoise_{t:N-1})$ in our policy to be a vector of scalar sigmoidal functions \[\varphi(\xi)=\frac{1- e^{-\xi}}{1+ e^{-\xi}} \] applied to each coordinate of the noise vector $\wnoise_t$. The decision variables $\offset_t, \gain_t$ and $\weight_t$ are computed at times $t = 0,\reachindex,2\reachindex,\cdots$, by solving optimization problems according to Theorems \ref{th:DF_rep_stbl} and \ref{th:DF_seq_stbl}. 
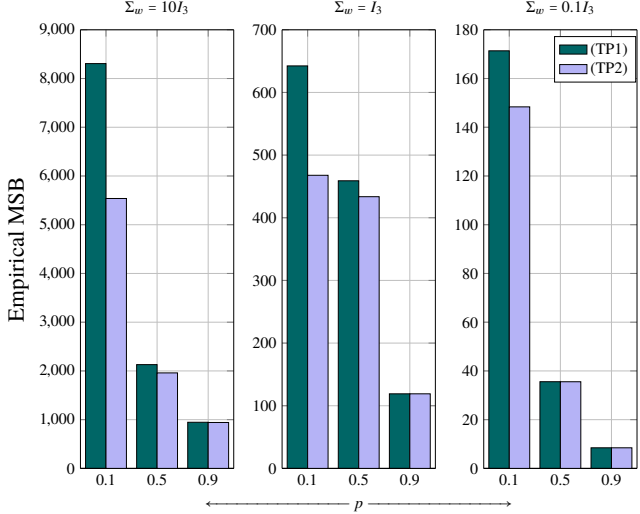
\begin{figure}[t]
		\centering
\begin{adjustbox}{width=\columnwidth}
%
%
\begin{tikzpicture}

\begin{axis}[%
width=1.328in,
height=3.803in,
at={(4.303in,0.513in)},
scale only axis,
separate axis lines,
every outer x axis line/.append style={black},
every x tick label/.append style={font=\color{black}},
xmin=0.5,
xmax=3.5,
xtick={1,2,3},
xticklabels={{0.1},{0.5},{0.9}},
xmajorgrids,
every outer y axis line/.append style={black},
every y tick label/.append style={font=\color{black}},
ymin=0,
ymax=180,
ymajorgrids,
axis background/.style={fill=white},
title = {$\Sigma_{\wnoise} = 0.1 I_3$},
legend style={at={(0.46,0.88)},anchor=south west,legend cell align=left,align=left,draw=black}
]
\addplot[ybar,bar width=0.4,draw=black,fill=mycolor1,area legend] plot table[row sep=crcr] {%
0.8	171.364228487467\\
1.8	35.5756258990014\\
2.8	8.45297986875135\\
};
\addlegendentry{\ref{a:seq}};

\addplot [color=black,solid,forget plot]
  table[row sep=crcr]{%
0.5	0\\
3.5	0\\
};
\addplot[ybar,bar width=0.4,draw=black,fill=mycolor2,area legend] plot table[row sep=crcr] {%
1.2	148.365805271095\\
2.2	35.5612087806662\\
3.2	8.45303904087852\\
};
\addlegendentry{\ref{a:repetitive}};

\end{axis}

\begin{axis}[%
width=1.328in,
height=3.803in,
at={(2.556in,0.513in)},
scale only axis,
separate axis lines,
every outer x axis line/.append style={black},
every x tick label/.append style={font=\color{black}},
xmin=0.5,
xmax=3.5,
xtick={1,2,3},
xticklabels={{0.1},{0.5},{0.9}},
xlabel={$\xleftarrow{\hspace*{3cm}} p \xrightarrow{\hspace*{3cm}}$},
xmajorgrids,
every outer y axis line/.append style={black},
every y tick label/.append style={font=\color{black}},
ymin=0,
ymax=700,
ymajorgrids,
axis background/.style={fill=white},
title = {$\Sigma_{\wnoise} = I_3$}
]
\addplot[ybar,bar width=0.4,draw=black,fill=mycolor1,area legend] plot table[row sep=crcr] {%
0.8	642.354960788564\\
1.8	459.0102518039\\
2.8	119.014798158919\\
};
\addplot [color=black,solid,forget plot]
  table[row sep=crcr]{%
0.5	0\\
3.5	0\\
};
\addplot[ybar,bar width=0.4,draw=black,fill=mycolor2,area legend] plot table[row sep=crcr] {%
1.2	467.790745394756\\
2.2	433.627245462715\\
3.2	119.029838576518\\
};
\end{axis}

\begin{axis}[%
width=1.328in,
height=3.803in,
at={(0.809in,0.513in)},
scale only axis,
separate axis lines,
every outer x axis line/.append style={black},
every x tick label/.append style={font=\color{black}},
xmin=0.5,
xmax=3.5,
xtick={1,2,3},
xticklabels={{0.1},{0.5},{0.9}},
xmajorgrids,
every outer y axis line/.append style={black},
every y tick label/.append style={font=\color{black}},
ymin=0,
ymax=9000,
ylabel={\large{Empirical MSB}},
ymajorgrids,
axis background/.style={fill=white},
title = {$\Sigma_{\wnoise} = 10 I_3$}
]
\addplot[ybar,bar width=0.4,draw=black,fill=mycolor1,area legend] plot table[row sep=crcr] {%
0.8	8306.22328184664\\
1.8	2129.34664675826\\
2.8	946.734746854056\\
};
\addplot [color=black,solid,forget plot]
  table[row sep=crcr]{%
0.5	0\\
3.5	0\\
};
\addplot[ybar,bar width=0.4,draw=black,fill=mycolor2,area legend] plot table[row sep=crcr] {%
1.2	5538.35010179029\\
2.2	1959.45333116425\\
3.2	942.207310112824\\
};
\end{axis}

\end{tikzpicture}%
\end{adjustbox}
\caption{Empirical mean square bound with $\mu = 1000$}
\label{Fig:msb_sparse_stbl}
\end{figure}
Our observations from the simulations are listed below. All quantities reported below correspond to averages over $100$ sample paths. 
\begin{enumerate}[leftmargin = *]
\item The empirical mean square bound (MSB) for \ref{a:seq} and \ref{a:repetitive} is plotted in Fig. \ref{Fig:msb_sparse_stbl} for $\mu = 1000$. The difference in empirical MSB between \ref{a:repetitive} and \ref{a:seq} increases with increase in variance and dropouts. When the successful transmission probability approaches 1 (less dropouts), the empirical MSB for \ref{a:repetitive} and \ref{a:seq} becomes equal.
\item The empirical average actuator energy is plotted in Fig.\ \ref{Fig:EnergyDF} for $\mu = 1000$. When the process noise variance is small, the actuator energy decreases with increase in $p$. For large value of the variance of additive process noise, it increases with increase in $p$. 
\item The computed control is sparse in the sense that there are about $7-19 \% $ time instants when null control is computed. See Fig.\ \ref{Fig:sparsity}. We observe that the percentage sparsity is large when successful transmission probability is less. When there are more dropouts due to interference and fading, the controller generates sparser control to reduce the impact of interference and fading. 
\end{enumerate}

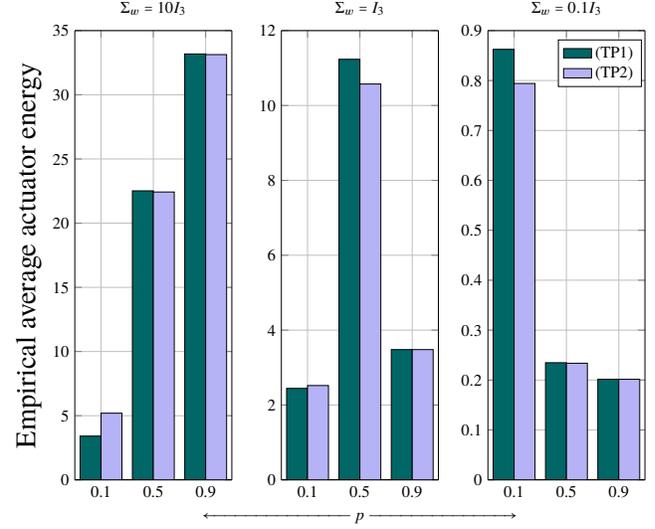
\begin{figure}[t]
		\centering
\begin{adjustbox}{width=\columnwidth}
%
%
\begin{tikzpicture}

\begin{axis}[%
width=1.328in,
height=3.803in,
at={(4.303in,0.513in)},
scale only axis,
separate axis lines,
every outer x axis line/.append style={black},
every x tick label/.append style={font=\color{black}},
xmin=0.5,
xmax=3.5,
xtick={1,2,3},
xticklabels={{0.1},{0.5},{0.9}},
xmajorgrids,
every outer y axis line/.append style={black},
every y tick label/.append style={font=\color{black}},
ymin=0,
ymax=0.9,
ymajorgrids,
axis background/.style={fill=white},
title={$\Sigma_{\wnoise} = 0.1 I_3$},
legend style={legend cell align=left,align=left,draw=black}
]
\addplot[ybar,bar width=0.4,draw=black,fill=mycolor1,area legend] plot table[row sep=crcr] {%
0.8	0.86249547336672\\
1.8	0.234729992360678\\
2.8	0.201497175246539\\
};
\addlegendentry{\ref{a:seq}};

\addplot [color=black,solid,forget plot]
  table[row sep=crcr]{%
0.5	0\\
3.5	0\\
};
\addplot[ybar,bar width=0.4,draw=black,fill=mycolor2,area legend] plot table[row sep=crcr] {%
1.2	0.794128527846899\\
2.2	0.233686012023938\\
3.2	0.201523538170634\\
};
\addlegendentry{\ref{a:repetitive}};

\end{axis}

\begin{axis}[%
width=1.328in,
height=3.803in,
at={(2.556in,0.513in)},
scale only axis,
separate axis lines,
every outer x axis line/.append style={black},
every x tick label/.append style={font=\color{black}},
xmin=0.5,
xmax=3.5,
xtick={1,2,3},
xticklabels={{0.1},{0.5},{0.9}},
xlabel={$ \xleftarrow{\hspace*{3cm}} p \xrightarrow{\hspace*{3cm}}$},
xmajorgrids,
every outer y axis line/.append style={black},
every y tick label/.append style={font=\color{black}},
ymin=0,
ymax=12,
ymajorgrids,
axis background/.style={fill=white},
title={$\Sigma_{\wnoise} = I_3$}
]
\addplot[ybar,bar width=0.4,draw=black,fill=mycolor1,area legend] plot table[row sep=crcr] {%
0.8	2.44461068356089\\
1.8	11.2341161058141\\
2.8	3.47968831953357\\
};
\addplot [color=black,solid,forget plot]
  table[row sep=crcr]{%
0.5	0\\
3.5	0\\
};
\addplot[ybar,bar width=0.4,draw=black,fill=mycolor2,area legend] plot table[row sep=crcr] {%
1.2	2.52003363586822\\
2.2	10.5737002944865\\
3.2	3.48008374002134\\
};
\end{axis}

\begin{axis}[%
width=1.328in,
height=3.803in,
at={(0.809in,0.513in)},
scale only axis,
separate axis lines,
every outer x axis line/.append style={black},
every x tick label/.append style={font=\color{black}},
xmin=0.5,
xmax=3.5,
xtick={1,2,3},
xticklabels={{0.1},{0.5},{0.9}},
xmajorgrids,
every outer y axis line/.append style={black},
every y tick label/.append style={font=\color{black}},
ymin=0,
ymax=35,
ylabel={\Large{Empirical average actuator energy}},
ymajorgrids,
axis background/.style={fill=white},
title={$\Sigma_{\wnoise} = 10 I_3$}
]
\addplot[ybar,bar width=0.4,draw=black,fill=mycolor1,area legend] plot table[row sep=crcr] {%
0.8	3.41259079880606\\
1.8	22.5197009069728\\
2.8	33.1793198413388\\
};
\addplot [color=black,solid,forget plot]
  table[row sep=crcr]{%
0.5	0\\
3.5	0\\
};
\addplot[ybar,bar width=0.4,draw=black,fill=mycolor2,area legend] plot table[row sep=crcr] {%
1.2	5.19632078152468\\
2.2	22.4280657661475\\
3.2	33.1320280561665\\
};
\end{axis}
\end{tikzpicture}%
\end{adjustbox}
\caption{Empirical average actuator energy for \ref{a:seq} and \ref{a:repetitive} with $\mu = 1000$}
\label{Fig:EnergyDF}
\end{figure}
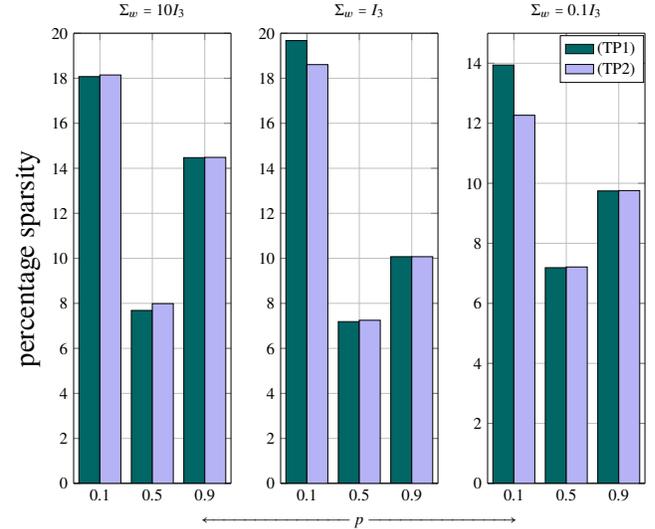
\begin{figure}[t]
		\centering
\begin{adjustbox}{width=\columnwidth}
%
%
\begin{tikzpicture}

\begin{axis}[%
width=1.328in,
height=3.803in,
at={(4.303in,0.513in)},
scale only axis,
separate axis lines,
every outer x axis line/.append style={black},
every x tick label/.append style={font=\color{black}},
xmin=0.5,
xmax=3.5,
xtick={1,2,3},
xticklabels={{0.1},{0.5},{0.9}},
xmajorgrids,
every outer y axis line/.append style={black},
every y tick label/.append style={font=\color{black}},
ymin=0,
ymax=15,
ymajorgrids,
axis background/.style={fill=white},
title = {$\Sigma_{\wnoise} = 0.1 I_3$},
legend style={at={(0.46,0.89)},anchor=south west,legend cell align=left,align=left,draw=black}
]
\addplot[ybar,bar width=0.4,draw=black,fill=mycolor1,area legend] plot table[row sep=crcr] {%
0.8	13.9393939393939\\
1.8	7.18787878787879\\
2.8	9.74545454545454\\
};
\addlegendentry{\ref{a:seq}};

\addplot [color=black,solid,forget plot]
  table[row sep=crcr]{%
0.5	0\\
3.5	0\\
};
\addplot[ybar,bar width=0.4,draw=black,fill=mycolor2,area legend] plot table[row sep=crcr] {%
1.2	12.2666666666667\\
2.2	7.21212121212121\\
3.2	9.75757575757576\\
};
\addlegendentry{\ref{a:repetitive}};

\end{axis}

\begin{axis}[%
width=1.328in,
height=3.803in,
at={(2.556in,0.513in)},
scale only axis,
separate axis lines,
every outer x axis line/.append style={black},
every x tick label/.append style={font=\color{black}},
xmin=0.5,
xmax=3.5,
xtick={1,2,3},
xticklabels={{0.1},{0.5},{0.9}},
xlabel = {$\xleftarrow{\hspace*{3cm}} p \xrightarrow{\hspace*{3cm}}$},
xmajorgrids,
every outer y axis line/.append style={black},
every y tick label/.append style={font=\color{black}},
ymin=0,
ymax=20,
ymajorgrids,
axis background/.style={fill=white},
title = {$\Sigma_{\wnoise} = I_3$}
]
\addplot[ybar,bar width=0.4,draw=black,fill=mycolor1,area legend] plot table[row sep=crcr] {%
0.8	19.6727272727273\\
1.8	7.18787878787879\\
2.8	10.0727272727273\\
};
\addplot [color=black,solid,forget plot]
  table[row sep=crcr]{%
0.5	0\\
3.5	0\\
};
\addplot[ybar,bar width=0.4,draw=black,fill=mycolor2,area legend] plot table[row sep=crcr] {%
1.2	18.6060606060606\\
2.2	7.24848484848485\\
3.2	10.0727272727273\\
};
\end{axis}

\begin{axis}[%
width=1.328in,
height=3.803in,
at={(0.809in,0.513in)},
scale only axis,
separate axis lines,
every outer x axis line/.append style={black},
every x tick label/.append style={font=\color{black}},
xmin=0.5,
xmax=3.5,
xtick={1,2,3},
xticklabels={{0.1},{0.5},{0.9}},
xmajorgrids,
every outer y axis line/.append style={black},
every y tick label/.append style={font=\color{black}},
ymin=0,
ymax=20,
ylabel={\Large{percentage sparsity}},
ymajorgrids,
axis background/.style={fill=white},
title = {$\Sigma_{\wnoise} = 10 I_3$}
]
\addplot[ybar,bar width=0.4,draw=black,fill=mycolor1,area legend] plot table[row sep=crcr] {%
0.8	18.0727272727273\\
1.8	7.68484848484849\\
2.8	14.4727272727273\\
};
\addplot [color=black,solid,forget plot]
  table[row sep=crcr]{%
0.5	0\\
3.5	0\\
};
\addplot[ybar,bar width=0.4,draw=black,fill=mycolor2,area legend] plot table[row sep=crcr] {%
1.2	18.1454545454545\\
2.2	7.98787878787879\\
3.2	14.4848484848485\\
};
\end{axis}

\end{tikzpicture}%
\end{adjustbox}
\caption{Percentage of time instants when null control is computed with $\mu = 1000$}
\label{Fig:sparsity}
\end{figure}

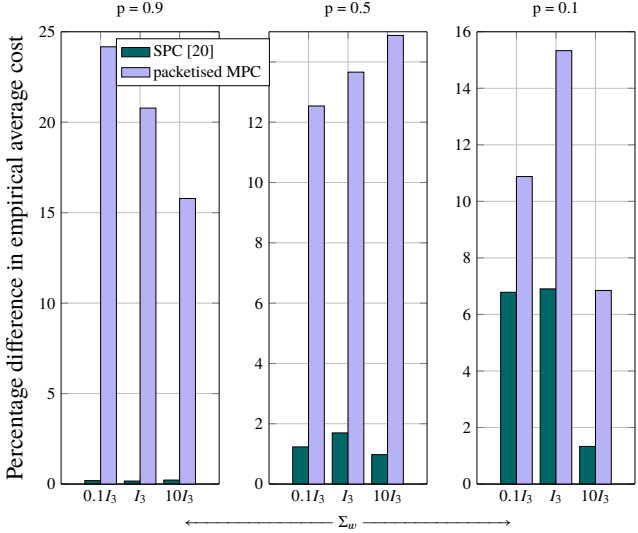
\begin{figure}[t]
		\centering
\begin{adjustbox}{width=\columnwidth}
%
%
\begin{tikzpicture}

\begin{axis}[%
width=1.328in,
height=3.803in,
at={(4.303in,0.513in)},
scale only axis,
separate axis lines,
every outer x axis line/.append style={black},
every x tick label/.append style={font=\color{black}},
xmin=0,
xmax=4,
xtick={1,2,3},
xticklabels={{$0.1 I_3$},{$I_3$},{$10 I_3$}},
xmajorgrids,
every outer y axis line/.append style={black},
every y tick label/.append style={font=\color{black}},
ymin=0,
ymax=16,
ymajorgrids,
axis background/.style={fill=white},
title={p = 0.1}
]
\addplot[ybar,bar width=0.4,draw=black,fill=mycolor1,area legend] plot table[row sep=crcr] {%
0.8	6.78377347072376\\
1.8	6.90388469662762\\
2.8	1.33011742514298\\
};
\addplot [color=black,solid,forget plot]
  table[row sep=crcr]{%
0	0\\
4	0\\
};
\addplot[ybar,bar width=0.4,draw=black,fill=mycolor2,area legend] plot table[row sep=crcr] {%
1.2	10.8774030290648\\
2.2	15.3310866845463\\
3.2	6.84872108523839\\
};
\end{axis}

\begin{axis}[%
width=1.328in,
height=3.803in,
at={(2.556in,0.513in)},
scale only axis,
separate axis lines,
every outer x axis line/.append style={black},
every x tick label/.append style={font=\color{black}},
xmin=0,
xmax=4,
xtick={1,2,3},
xticklabels={{$0.1 I_3$},{$I_3$},{$10 I_3$}},
xmajorgrids,
xlabel = {$\xleftarrow{\hspace*{3cm}} \Sigma_{\wnoise} \xrightarrow{\hspace*{3cm}}$},
every outer y axis line/.append style={black},
every y tick label/.append style={font=\color{black}},
ymin=0,
ymax=15,
ymajorgrids,
axis background/.style={fill=white},
title={p = 0.5}
]
\addplot[ybar,bar width=0.4,draw=black,fill=mycolor1,area legend] plot table[row sep=crcr] {%
0.8	1.23092839958871\\
1.8	1.69600855008283\\
2.8	0.975396345360887\\
};
\addplot [color=black,solid,forget plot]
  table[row sep=crcr]{%
0	0\\
4	0\\
};
\addplot[ybar,bar width=0.4,draw=black,fill=mycolor2,area legend] plot table[row sep=crcr] {%
1.2	12.5380830538921\\
2.2	13.6609903699536\\
3.2	14.8786259467396\\
};
\end{axis}

\begin{axis}[%
width=1.328in,
height=3.803in,
at={(0.809in,0.513in)},
scale only axis,
separate axis lines,
every outer x axis line/.append style={black},
every x tick label/.append style={font=\color{black}},
xmin=0,
xmax=4,
xtick={1,2,3},
xticklabels={{$0.1 I_3$},{$I_3$},{$10 I_3$}},
xmajorgrids,
every outer y axis line/.append style={black},
every y tick label/.append style={font=\color{black}},
ymin=0,
ymax=25,
ylabel={\large{Percentage difference in empirical average cost}},
ymajorgrids,
axis background/.style={fill=white},
title={p = 0.9},
legend style={at={(0.348,0.88)},anchor=south west,legend cell align=left,align=left,draw=black}
]
\addplot[ybar,bar width=0.4,draw=black,fill=mycolor1,area legend] plot table[row sep=crcr] {%
0.8	0.19047770001129\\
1.8	0.164364452794474\\
2.8	0.219406667072399\\
};
\addlegendentry{SPC \cite{ref:PDQ-15}};

\addplot [color=black,solid,forget plot]
  table[row sep=crcr]{%
0	0\\
4	0\\
};
\addplot[ybar,bar width=0.4,draw=black,fill=mycolor2,area legend] plot table[row sep=crcr] {%
1.2	24.1687874370294\\
2.2	20.7803732979079\\
3.2	15.7848132006999\\
};
\addlegendentry{packetised MPC};

\end{axis}
\end{tikzpicture}%
\end{adjustbox}
\caption{Percentage difference in empirical average cost for \ref{a:repetitive} with respect to proposed approach}
\label{Fig:CostRepDF}
\end{figure}

\begin{figure}[t]
		\centering
\begin{adjustbox}{width=\columnwidth}
%
%
\begin{tikzpicture}

\begin{axis}[%
width=1.328in,
height=3.803in,
at={(4.303in,0.513in)},
scale only axis,
separate axis lines,
every outer x axis line/.append style={black},
every x tick label/.append style={font=\color{black}},
xmin=0,
xmax=4,
xtick={1,2,3},
xticklabels={{$0.1 I_3$},{$I_3$},{$10 I_3$}},
xmajorgrids,
every outer y axis line/.append style={black},
every y tick label/.append style={font=\color{black}},
ymin=-4,
ymax=12,
ymajorgrids,
axis background/.style={fill=white},
title={p = 0.1}
]
\addplot[ybar,bar width=0.4,draw=black,fill=mycolor1,area legend] plot table[row sep=crcr] {%
0.8	-2.40203067007098\\
1.8	4.55218707608812\\
2.8	1.30821790629594\\
};
\addplot [color=black,solid,forget plot]
  table[row sep=crcr]{%
0	0\\
4	0\\
};
\addplot[ybar,bar width=0.4,draw=black,fill=mycolor2,area legend] plot table[row sep=crcr] {%
1.2	1.48743625860795\\
2.2	10.6683453577748\\
3.2	2.73040516887357\\
};
\end{axis}

\begin{axis}[%
width=1.328in,
height=3.803in,
at={(2.556in,0.513in)},
scale only axis,
separate axis lines,
every outer x axis line/.append style={black},
every x tick label/.append style={font=\color{black}},
xmin=0,
xmax=4,
xtick={1,2,3},
xticklabels={{$0.1 I_3$},{$I_3$},{$10 I_3$}},
xmajorgrids,
every outer y axis line/.append style={black},
every y tick label/.append style={font=\color{black}},
ymin=0,
ymax=18,
ymajorgrids,
xlabel = {$\xleftarrow{\hspace*{3cm}} \Sigma_{\wnoise} \xrightarrow{\hspace*{3cm}}$},
axis background/.style={fill=white},
title={p = 0.5}
]
\addplot[ybar,bar width=0.4,draw=black,fill=mycolor1,area legend] plot table[row sep=crcr] {%
0.8	4.19518048084799\\
1.8	4.45607488101719\\
2.8	1.59930010016336\\
};
\addplot [color=black,solid,forget plot]
  table[row sep=crcr]{%
0	0\\
4	0\\
};
\addplot[ybar,bar width=0.4,draw=black,fill=mycolor2,area legend] plot table[row sep=crcr] {%
1.2	17.5160096787328\\
2.2	14.65443587657\\
3.2	17.085610686637\\
};
\end{axis}

\begin{axis}[%
width=1.328in,
height=3.803in,
at={(0.809in,0.513in)},
scale only axis,
separate axis lines,
every outer x axis line/.append style={black},
every x tick label/.append style={font=\color{black}},
xmin=0,
xmax=4,
xtick={1,2,3},
xticklabels={{$0.1 I_3$},{$I_3$},{$10 I_3$}},
xmajorgrids,
every outer y axis line/.append style={black},
every y tick label/.append style={font=\color{black}},
ymin=0,
ymax=30,
ylabel={\large{Percentage difference in empirical average cost}},
ymajorgrids,
axis background/.style={fill=white},
title={p = 0.9},
legend style={at={(0.13,0.88)},anchor=south west,legend cell align=left,align=left,draw=black}
]
\addplot[ybar,bar width=0.4,draw=black,fill=mycolor1,area legend] plot table[row sep=crcr] {%
0.8	0.247911522456703\\
1.8	0.431505326334589\\
2.8	0.543258208586762\\
};
\addlegendentry{SPC \cite{ref:PDQ-15}};

\addplot [color=black,solid,forget plot]
  table[row sep=crcr]{%
0	0\\
4	0\\
};
\addplot[ybar,bar width=0.4,draw=black,fill=mycolor2,area legend] plot table[row sep=crcr] {%
1.2	25.7411785137798\\
2.2	21.7871396065871\\
3.2	17.4667739607953\\
};
\addlegendentry{MPC};

\end{axis}
\end{tikzpicture}%
\end{adjustbox}
\caption{Percentage difference in empirical average cost for \ref{a:seq} with respect to proposed approach}
\label{Fig:CostSeqDF}
\end{figure}
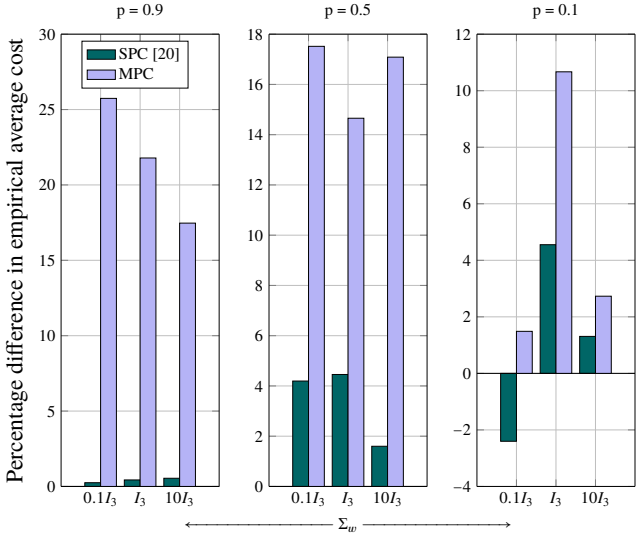
\par In CE approaches (PPC, MPC and packetised MPC), the control sequence is generated by considering nominal plant model and applied to the stochastic system \eqref{e:example}. In MPC only present control value is transmitted through channel, similar to \ref{a:seq}. But in packetised MPC, future control values are also transmitted until first successful transmission, similar to \ref{a:repetitive}. In SPC, MPC and packetized MPC, the optimization is carried out after each $N_r$ time steps, but in PPC, the optimization is carried out after each time step. In \cite{ref:PDQ-15} SPC and PPC \cite{quevedo-12} are compared and it is shown that SPC \cite{ref:PDQ-15} outperforms PPC \cite{quevedo-12}. Here we compare our approach with SPC \cite{ref:PDQ-15} and with CE approaches by forcing $\mu = 0$ in Theorems \ref{th:DF_rep_stbl} and \ref{th:DF_seq_stbl}.  We have following observations: 
\begin{enumerate}[leftmargin  = *, start=4]
\item The present approach and SPC \cite{ref:PDQ-15} outperform CE approaches in all cases. See Fig.\ \ref{Fig:CostRepDF} and \ref{Fig:CostSeqDF}.
\item Our approach performs better than SPC \cite{ref:PDQ-15} in all considered cases except in one case when successful transmission probability is very low $(p = 0.1)$, variance of additive process noise is small $\Sigma_{\wnoise}=0.1 I_3$ and sequential transmission protocol \ref{a:seq} is used. See Fig.\ \ref{Fig:CostSeqDF}.
\item The percentage difference in empirical average cost between our approach and SPC \cite{ref:PDQ-15} is large when successful transmission probability is moderate $p=0.5$ and sequential transmission \ref{a:seq} is used.
\item For repetitive transmission protocol \ref{a:repetitive}, the percentage difference in empirical average cost between our approach and SPC \cite{ref:PDQ-15} is large when the successful transmission probability $p$ is very small.    
\end{enumerate}
 
\par
\begin{figure}[t]
\centering
\begin{adjustbox}{width=0.8\columnwidth}
\input{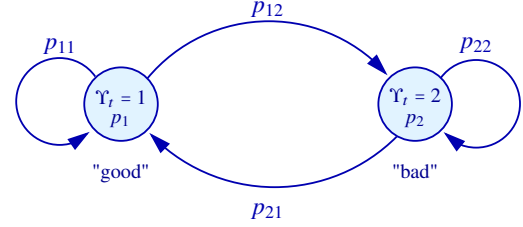}
\end{adjustbox}
\caption{Transmission dropout model with a binary network state
    $(\Upsilon_t)_{t\in \{1,2\}}$: when $\Upsilon_t=1$ the channel is reliable with high successful transmission probabilities; $\Upsilon_t=2$ refers to a situation where the channel is
  unreliable and transmissions are more likely to be dropped.}
  \label{Fig:Gilbert-Eliot}
\end{figure}

\begin{figure}[t]
		\centering
\begin{adjustbox}{width=\columnwidth}
%
%
\begin{tikzpicture}

\begin{axis}[%
width=4.822in,
height=3.803in,
at={(0.809in,0.513in)},
scale only axis,
separate axis lines,
every outer x axis line/.append style={black},
every x tick label/.append style={font=\color{black}},
xmin=0.5,
xmax=1,
xtick={0.5,0.6,0.7,0.8,0.9,1},
xlabel={$\text{p}_\text{1}$},
xmajorgrids,
every outer y axis line/.append style={black},
every y tick label/.append style={font=\color{black}},
ymin=65,
ymax=110,
ylabel={\Large{Empirical average cost}},
ymajorgrids,
axis background/.style={fill=white},
legend style={at={(0.605,0.83)},anchor=south west,legend cell align=left,align=left,draw=black}
]
\addplot [color=red,solid,line width=2.0pt]
  table[row sep=crcr]{%
0.5	95.7979915959223\\
0.6	87.3204638223328\\
0.7	80.2212253960004\\
0.8	74.7886157648109\\
0.9	66.8988306313955\\
1	65.0563923086126\\
};
\addlegendentry{SPC \ref{a:repetitive}};

\addplot [color=blue,solid,line width=2.0pt]
  table[row sep=crcr]{%
0.5	106.311373097075\\
0.6	96.805388397306\\
0.7	89.7276092533935\\
0.8	84.484575079412\\
0.9	76.8590075910678\\
1	75.1623136934371\\
};
\addlegendentry{Packetised MPC};

\end{axis}
\end{tikzpicture}%
\end{adjustbox}
\caption{Empirical average cost for \ref{a:repetitive} with correlated channel noise}
\label{Fig:repCorrelated}
\end{figure}
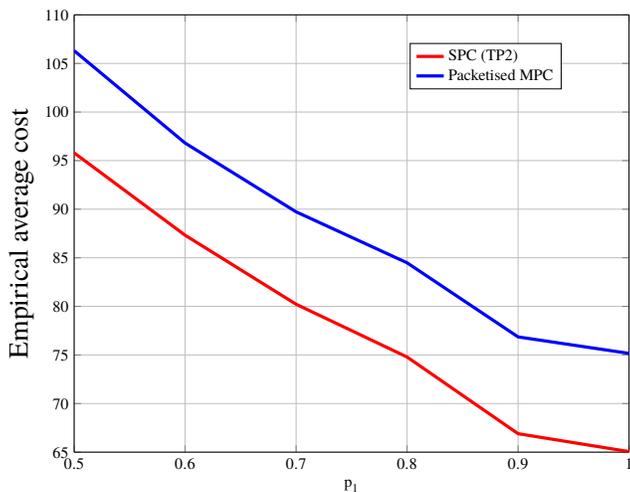

 We have tested our approach on correlated channel noise also by considering the Gilbert-Eliot channel model \cite{gilbert} as given in Fig.\ \ref{Fig:Gilbert-Eliot}. We have fixed successful transmission probability for the bad channel $p_2 = 0$, the transmission probabilities from good to bad channel $p_{12} = 0.2$, from bad to good channel $p_{21} = 0.9$. We have simulated for the successful transmission probability of good channel $p_1 \in \{ 0.5, 0.6, 0.7, 0.8, 0.9, 1 \}$ and variance of additive process noise $\Sigma_{\wnoise} = 5 I_3$ with given dynamics of the plant \eqref{e:example}. The plots for \ref{a:seq} and \ref{a:repetitive} are shown in Fig.\ \ref{Fig:seqCorrelated} and \ref{Fig:repCorrelated}, respectively. The transmission protocol \ref{a:seq} is compared with MPC in which no buffer is used at the actuator end. The transmission protocol \ref{a:repetitive} is compared with packetised MPC in which a buffer of the size of optimization horizon is used at the actuator end. For the purpose of fair comparison the recalculation interval is fixed $N_r = 3$ in all cases. We can observe that our present approach performs better than older approaches.   

\begin{figure}[t]
		\centering
\begin{adjustbox}{width=\columnwidth}
%
%
\begin{tikzpicture}

\begin{axis}[%
width=4.822in,
height=3.803in,
at={(0.809in,0.513in)},
scale only axis,
separate axis lines,
every outer x axis line/.append style={black},
every x tick label/.append style={font=\color{black}},
xmin=0.5,
xmax=1,
xtick={0.5,0.6,0.7,0.8,0.9,1},
xlabel={$\text{p}_\text{1}$},
xmajorgrids,
every outer y axis line/.append style={black},
every y tick label/.append style={font=\color{black}},
ymin=60,
ymax=160,
ylabel={\Large{Empirical average cost}},
ymajorgrids,
axis background/.style={fill=white},
legend style={at={(0.658,0.828)},anchor=south west,legend cell align=left,align=left,draw=black}
]
\addplot [color=red,solid,line width=2.0pt]
  table[row sep=crcr]{%
0.5	138.750257495848\\
0.6	111.149506730119\\
0.7	95.6535759511996\\
0.8	85.8196515949619\\
0.9	75.7987008834804\\
1	67.5186629850579\\
};
\addlegendentry{SPC \ref{a:seq}};

\addplot [color=blue,solid,line width=2.0pt]
  table[row sep=crcr]{%
0.5	158.828202884916\\
0.6	128.775587142288\\
0.7	109.21692663785\\
0.8	99.4057594872696\\
0.9	91.9182746878713\\
1	80.3336421651774\\
};
\addlegendentry{MPC};

\end{axis}
\end{tikzpicture}%
\end{adjustbox}
\caption{Empirical average cost for \ref{a:seq} with correlated channel noise}
\label{Fig:seqCorrelated}
\end{figure}
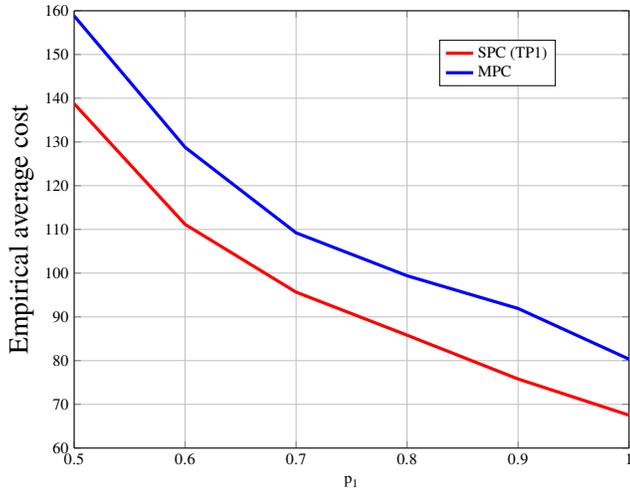

%
%

\section{Epilogue}\label{s:epilogue}
We have demonstrated that the proposed feedback policy leads to a convex quadratic program and the introduction of sparsity helps in communication-control co-design. We have shown by numerical experiments how the matrices involved in the construction of the objective function are computed offline. The current approach is restricted to i.i.d. packet dropouts and assumes that complete state information is available. The extensions of the ideas presented here may include multi-channel systems \cite{koegel2013distributed,ljevsnjanin2014packetized}, self triggered or event triggered operations \cite{araujo2014system,donkers2012output,Nagahara16TAC,quevedo2014stochastic,prabhat2015event}. Our approach can also be extended for the case of more general multiplicative noise. 


\bibliographystyle{IEEEtran}        
\bibliography{refs}           
\appendix

\section{Proofs}\label{s:proofs}         
	This appendix provides the proofs of our theorems.
	 Let us  define the component-wise saturation function  \(\R^{d_o}\ni z\longmapsto \sat_{r, \zeta}^\infty(z)\in\R^{d_o}\) to be 
				\[
					\bigl(\sat_{r, \zeta}^\infty(z)\bigr)_{i} = \begin{cases}
						z_i \zeta/r	& \text{if \(\abs{z_i} \le r\),}\\
						\zeta	 	& \text{if \(z_i > r\), and }\\
						-\zeta		& \text{otherwise,}
					\end{cases}
				\]
				for each \(i = 1, \ldots, d_o\).
\begin{lem}
\label{t:bhatia}
The matrix $A$ is positive semidefinite if and only if $\calA = \pmat{A & \cdots & A\\ \vdots & \ddots & \vdots \\ A & \cdots & A}$ is positive semidefinite.
\end{lem}	
\begin{pf}
The proof is in same line of arguments as in \cite[Lemma 1.3.6]{bhatia2009positive}. We are providing the proof for the completeness. Let $A$ be positive semidefinite. Then, we have \[ \calA = \pmat{A^{1/2} & 0 & \cdots & 0 \\ \vdots & \vdots & \ddots & \vdots \\ A^{1/2} & 0 & \cdots & 0 } \pmat{A^{1/2} & \cdots & A^{1/2} \\ 0 & \cdots & 0 \\ \vdots & \ddots & \vdots \\ 0 & \cdots & 0} ,\]
which in turn implies $\calA$ is also positive semi-definite. Now, if we consider $\calA$ to be positive semi-definite, we get $\calA = C\transp C$ for some unique positive semi-definite matrix $C$. By comparing, the diagonal block of left and right hand side, we can conclude that $A$ is positive semi-definite.  
\end{pf}					
\begin{lem} \label{l:objective}
Consider the system \eqref{e:system}, for every $t=0, \reachindex, 2\reachindex, \cdots,$ the problem \eqref{e:opt problem} under policy \eqref{e:policy}, transmission protocol \ref{a:repetitive} and control set \eqref{e:controlset} is convex quadratic with respect to the decision variable $\Xi_t$ and $\gain_t$. The objective function \eqref{e:opt problem sparse} is given by \eqref{e:programsingle}. 
\end{lem}	
\begin{pf}
The objective function of \eqref{e:opt problem} is given by
\begin{align*}
&\EE_{\st_{t}}\bigl[ \inprod{\st_{t:N+1}}{\calQ \st_{t:N+1}} + \inprod{\control^a_{t:N}}{\calR\control^a_{t:N}}\bigr]\\
&=\EE_{\st_{t}}\Bigl[\inprod{\calA \st_t + \calB \control^a_{t:N} + \calD\wnoise_{t:N} }{ \calQ \Bigl( \calA \st_t  + \calB \control^a_{t:N} + \calD\wnoise_{t:N} \Bigr)} \\
 & \quad + \inprod{\control^a_{t:N}}{ \calR\control^a_{t:N} } \Bigr]\\
&=\EE_{\st_{t}}\Bigl[\inprod{\calA \st_t}{\calQ \calA \st_t} + \inprod{\calD\wnoise_{t:N}}{\calQ\calD\wnoise_{t:N} }+ 2\inprod{ \calA \st_t }{ \calQ  \calB \control^a_{t:N} } \\ 
 & \quad + 2\EE_{\st_{t}}\bigl[ (\calS\gain_t\ee(\wnoise_{t:N-1}))\transp\calB\transp\calQ\calD\wnoise_{t:N} \bigr] \\ 
 & \quad + \inprod{\control^a_{t:N}}{ (\calB\transp\calQ\calB+\calR )\control^a_{t:N} } \Bigr]\\
        &= \inprod{\calA \st_t}{\calQ \calA \st_t} + \trace(\calD \transp \calQ \calD \Sigma_{\wnoise}) + 2\inprod{ \calA \st_t }{\calQ  \calB \mu_{\calG}\offset_t} \\
& \quad+ 2\EE_{\st_{t}}\bigl[ (\calS\gain_t\ee(\wnoise_{t:N-1}))\transp\calB\transp\calQ\calD\wnoise_{t:N} \bigr] + \inprod{\offset_t}{\Sigma_{\calG}\offset_t} \\ 
& \quad + \EE_{\st_t}\bigl[ \inprod{\calS \gain_t \ee(\wnoise_{t:N-1})}{\alpha \calS \gain_t \ee(\wnoise_{t:N-1})} \bigr] \\
& \quad + \EE_{\st_{t}}\Bigl[2\inprod{ \calA \st_t }{ \calQ  \calB \left( \calS\weight_t\cnoise_{t:N-1} \right) }  + 2\inprod{\offset_t}{\calG \transp\alpha \calS \weight_t\cnoise_{t:N-1}} \\ 
 & \quad+\inprod{\left( \weight_t\cnoise_{t:N-1} \right)}{ \calS \transp\alpha  \calS\left(\weight_t\cnoise_{t:N-1} \right) } \Bigr] \\
        &= \inprod{\calA \st_t}{\calQ \calA \st_t} + \trace(\calD \transp \calQ \calD        \Sigma_{\wnoise}) + 2\inprod{ \calA \st_t }{\calQ  \calB \mu_{\calG}\offset_t} \\
& \quad + 2\trace(\gain_t\transp\mu_{\calS}\transp\calB\transp\calQ\calD\Sigma_{\ee}^{\prime}) + \inprod{\offset_t}{\Sigma_{\calG}\offset_t} + \trace(\gain_t \transp \Sigma_{\calS} \gain_t \Sigma_{\ee}) \\ 
& \quad + \EE_{\st_{t}}\Bigl[2\inprod{ \calA \st_t }{ \calQ  \calB \left( \calS\weight_t\cnoise_{t:N-1} \right) }  + 2\inprod{\offset_t}{\calG \transp\alpha \calS \weight_t\cnoise_{t:N-1}} \\ 
 & \quad+\inprod{\left( \weight_t\cnoise_{t:N-1} \right)}{ \calS \transp\alpha  \calS\left(\weight_t\cnoise_{t:N-1} \right) } \Bigr] \\
       &= \inprod{\calA \st_t}{\calQ \calA \st_t} + \trace(\calD \transp \calQ \calD \Sigma_{\wnoise}) + 2\inprod{ \calA \st_t }{\calQ  \calB \mu_{\calG}\offset_t}\\
& \quad + 2\trace(\gain_t\transp\mu_{\calS}\transp\calB\transp\calQ\calD\Sigma_{\ee}^{\prime})  + \inprod{\offset_t}{\Sigma_{\calG}\offset_t}  + \trace(\gain_t \transp \Sigma_{\calS} \gain_t \Sigma_{\ee}) \\ 
& \quad + \EE_{\st_{t}}\Bigg[ 2 \inprod{ \calA \st_t }{ \calQ  \calB \left( \calS\sum_{\ell=1}^{N-1}\cnoise_{t+\ell-1}\weight_t^{(:,\ell)} \right) } \\
& \quad + 2\inprod{\offset_t}{\calG \transp\alpha \calS \sum_{\ell=1}^{N-1}\cnoise_{t+\ell-1}\weight_t^{(:,\ell)}} \\ 
 & \quad+\inprod{\left( \sum_{n=1}^{N-1}\cnoise_{t+n-1}\weight_t^{(:,n)} \right)}{ \calS \transp\alpha  \calS\left(\sum_{\ell=1}^{N-1}\cnoise_{t+\ell-1}\weight_t^{(:,\ell)} \right) } \Bigg] \\
          &= \inprod{\calA \st_t}{\calQ \calA \st_t} + \trace(\calD \transp \calQ \calD \Sigma_{\wnoise}) + 2\inprod{ \calA \st_t }{\calQ  \calB \mu_{\calG}\offset_t} \\ 
 & \quad + 2\trace(\gain_t\transp\mu_{\calS}\transp\calB\transp\calQ\calD\Sigma_{\ee}^{\prime}) + \inprod{\offset_t}{\Sigma_{\calG}\offset_t}  + \trace(\gain_t \transp \Sigma_{\calS} \gain_t \Sigma_{\ee}) \\& \quad +2\sum_{\ell=1}^{N-1}\inprod{ \calA \st_t }{ \calQ  \calB \EE_{\st_{t}}\left[\cnoise_{t+\ell-1}\calS\right] \weight_t^{(:,\ell)} }   \\ 
& \quad + \sum_{\ell=1}^{N-1}  2\inprod{\offset_t}{\EE_{\st_{t}}\Bigl[\calG \transp\alpha \calS \cnoise_{t+\ell-1}\Bigr] \weight_t^{(:,\ell)}} \\ 
 & \quad+\sum_{n=1}^{N-1}\sum_{\ell=1}^{N-1}\inprod{\left( \weight_t^{(:,n)} \right)}{\EE_{\st_{t}}\Bigl[\cnoise_{t+n-1} \calS \transp\alpha  \calS \cnoise_{t+\ell-1} \Bigr]\weight_t^{(:,\ell)}  } \\
        &= \inprod{\calA \st_t}{\calQ \calA \st_t} + \trace(\calD \transp \calQ \calD \Sigma_{\wnoise}) + 2\inprod{ \calA \st_t }{\calQ  \calB \mu_{\calG}\offset_t} \\
& \quad + 2\trace(\gain_t\transp\mu_{\calS}\transp\calB\transp\calQ\calD\Sigma_{\ee}^{\prime}) + \inprod{\offset_t}{\Sigma_{\calG}\offset_t}  + \trace(\gain_t \transp \Sigma_{\calS} \gain_t \Sigma_{\ee}) \\ 
& \quad + 2\sum_{\ell=1}^{N-1}\inprod{ \calA \st_t }{ \calQ  \calB \EE_{\st_{t}}\left[\calS_{\ell} \right] \weight_t^{(:,\ell)} }   \\& \quad +  \sum_{\ell=1}^{N-1}  2\inprod{\offset_t}{\EE_{\st_{t}}\Bigl[\calG \transp\alpha \calS_{\ell} \Bigr] \weight_t^{(:,\ell)}} \\ 
 & \quad+\sum_{n=1}^{N-1}\sum_{\ell=1}^{N-1}\inprod{\left( \weight_t^{(:,n)} \right)}{\EE_{\st_{t}}\Bigl[\calS_{n} \transp\alpha  \calS_{\ell} \Bigr]\weight_t^{(:,\ell)}  } \\
        &= \inprod{\calA \st_t}{\calQ \calA \st_t}+ \trace(\calD \transp \calQ \calD \Sigma_{\wnoise}) + 2\inprod{ \calA \st_t }{\calQ  \calB \mu_{\calG}\offset_t} \\ 
& \quad + 2\trace(\gain_t\transp\mu_{\calS}\transp\calB\transp\calQ\calD\Sigma_{\ee}^{\prime}) + \inprod{\offset_t}{\Sigma_{\calG}\offset_t}  + \trace(\gain_t \transp \Sigma_{\calS} \gain_t \Sigma_{\ee}) \\
& \quad + 2\sum_{\ell=1}^{N-1}\inprod{ \calA \st_t }{ \calQ  \calB \mu_{\calS_{\ell}} \weight_t^{(:,\ell)} }  \\
& \quad  + 2\sum_{\ell=1}^{N-1} \inprod{\offset_t}{\Sigma_{\calS\calG_{\ell}} \weight_t^{(:,\ell)}} \\
& \quad +\sum_{n=1}^{N-1}\sum_{\ell=1}^{N-1}\inprod{ \weight_t^{(:,n)} }{\Sigma_{\calS_{ n \ell}}\weight_t^{(:,\ell)}  } \\
        &= \inprod{\calA \st_t}{\calQ \calA \st_t} + \trace(\calD \transp \calQ \calD \Sigma_{\wnoise}) + 2\inprod{ \calA \st_t }{\calQ  \calB \mu_{\calG}\offset_t} \\
& \quad + 2\trace(\gain_t\transp\mu_{\calS}\transp\calB\transp\calQ\calD\Sigma_{\ee}^{\prime}) +  \inprod{\offset_t}{\Sigma_{\calG}\offset_t}  + \trace(\gain_t \transp \Sigma_{\calS} \gain_t \Sigma_{\ee}) \\
& \quad +2\inprod{ \calA \st_t }{ \calQ  \calB \tilde{\mu}_{\calS_{\ell}} \tilde{\weight}_t}  + 2 \inprod{\offset_t}{\tilde{\Sigma}_{\calS\calG_{\ell}} \tilde{\weight}_t} + \inprod{ \tilde{\weight}_t }{\tilde{\Sigma}_{\calS_{ n \ell}}\tilde{\weight}_t  } .
\end{align*}	
By rearranging the terms, we get
\begin{align*}
&\EE_{\st_{t}}\bigl[ \inprod{\st_{t:N+1}}{\calQ \st_{t:N+1}} + \inprod{\control^a_{t:N}}{\calR\control^a_{t:N}}\bigr]\\
& = \inprod{\Xi_t}{ \pmat{\Sigma_{\calG} & \tilde{\Sigma}_{\calS\calG_{\ell}} \\ \tilde{\Sigma}_{\calS\calG_{\ell}}\transp & \tilde{\Sigma}_{\calS_{ n \ell}} } \Xi_t } + \inprod{\calA \st_t}{2 \pmat{\calQ \calB \mu_{\calG} & \calQ \calB \tilde{\mu}_{\calS_{\ell}} } \Xi_t } \\ & \quad +  \inprod{\calA \st_t}{\calQ \calA \st_t} + \trace(\calD \transp \calQ \calD \Sigma_{\wnoise}) \\
& \quad +  2\trace(\gain_t\transp\mu_{\calS}\transp\calB\transp\calQ\calD\Sigma_{\ee}^{\prime}) +\trace(\gain_t \transp \Sigma_{\calS} \gain_t \Sigma_{\ee})\\
& = \inprod{\Xi_t}{ \calL \Xi_t } + \inprod{\calM \calA \st_t}{  \Xi_t } +  \inprod{\calA \st_t}{\calQ \calA \st_t} \\
& \quad+ 2\trace(\gain_t\transp\mu_{\calS}\transp\calB\transp\calQ\calD\Sigma_{\ee}^{\prime}) + \trace(\calD \transp \calQ \calD \Sigma_{\wnoise}) \\
&\quad + \trace(\gain_t \transp \Sigma_{\calS} \gain_t \Sigma_{\ee}).
\end{align*}
Now, by including the regularization term \eqref{e:regularization term} we get the objective function of \eqref{e:programsingle}.
Let us define \\ $\mathfrak{C} \Let \blkdiag\bigl( \mathfrak{C}_0, \mathfrak{C}_1 , \cdots, \mathfrak{C}_{N-1} \bigr)$, $\tilde{\calS} \Let \blkdiag\bigl( \calG, \calS_1 , \cdots, \calS_{N-1} \bigr)$ and  \[\mathbb{P} = \pmat{\alpha & \cdots & \alpha \\ \vdots & \ddots & \vdots \\ \alpha & \cdots & \alpha
}, \]
where $\alpha = \calB\transp\calQ\calB+\calR$ is a symmetric positive definite matrix. The matrix $\mathbb{P}$ can be shown to be positive semi-definite using Lemma \ref{t:bhatia}. 
Now, we have
\begin{align*} 
 \calL &= \EE_{\st_t}[\mathfrak{C}\transp\tilde{\calS}\transp \mathbb{P} \tilde{\calS} \mathfrak{C} ] = \EE_{\st_t}[\mathfrak{C}\transp\tilde{\calS}\transp \mathbb{P}^{1/2} \mathbb{P}^{1/2} \tilde{\calS} \mathfrak{C} ] \\
&= \EE_{\st_t}[\mathfrak{C}\transp \mathbb{P}^{1/2} \tilde{\calS} \tilde{\calS} \transp \mathbb{P}^{1/2} \mathfrak{C} ] = \mathfrak{C}\transp \mathbb{P}^{1/2} \EE_{\st_t}[ \tilde{\calS} \tilde{\calS} \transp ] \mathbb{P}^{1/2} \mathfrak{C}, 
\end{align*} 
where $\EE[ \tilde{\calS} \tilde{\calS} \transp ] \teL \mu_{\tilde{\calS}}$ has zero off-diagonal entries and positive diagonal entries. Hence, $\mu_{\tilde{\calS}}$ is positive definite and $\calL$ is positive semi definite. The objective function in \eqref{e:programsingle} is quadratic in $\Xi_t$ and involves trace operator on $\gain_t$. The regularization term is obviously convex. Therefore the objective function in \eqref{e:programsingle} is convex quadratic. 
\end{pf}

\begin{pf}[Proof of Theorem \ref{th:DF_rep}]
The proof of convexity and the formulation of the objective are given in Lemma \ref{l:objective}. The proposed control policy \eqref{e:policy} satisfies hard input constraint \eqref{e:controlset} as long as the following condition is satisfied: \[ \norm{\offset_t + \weight_t \cnoise_{t:N-1} + \gain_t \ee(\wnoise_{t:N-1})}_{\infty} \leq U_{\max} \]  for all $\cnoise_{t:N-1} \in \{ 0,1\}^N$. This is equivalent to the condition 
\[
 \max_{\substack{\cnoise_{t:N-1} \in \{ 0,1\}^{N-1} \\ \norm{\ee(\wnoise_{t:N-1})}_{\infty} \leq  \varphi_{\max} }} \abs{ \offset_t^{(i)} + \weight_t^{(i,:)}\cnoise_{t:N-1} + \gain_t^{(i,:)}\ee(\wnoise_{t:N-1}) } \leq U_{\max} 
\]
for all $i = 1,\cdots, mN $. Let $g_t \Let \cnoise_{t:N-1} - \dfrac{1}{2} \ones_{N} $ and $h_i \Let \offset_t^{(i)} + \dfrac{1}{2} \weight_t^{(i,:)}\ones_{N} $. Then preceding inequality is equivalent to
\begin{equation*}
\left \{
\begin{aligned}
\max_{\substack{g_t \in \{ -1/2,1/2\}^{N-1} \\ \norm{\ee(\wnoise_{t:N-1})}_{\infty} \leq  \varphi_{\max} }}  h_i + \weight_t^{(i,:)}g_i + \gain_t^{(i,:)}\ee(\wnoise_{t:N-1})  & \leq U_{\max}  \\ 
 \min_{\substack{g_t \in \{ -1/2, 1/2 \}^{N-1} \\ \norm{\ee(\wnoise_{t:N-1})}_{\infty}  \leq  \varphi_{\max} }}  h_i + \weight_t^{(i,:)}g_i + \gain_t^{(i,:)}\ee(\wnoise_{t:N-1}) & \geq - U_{\max} 
 \end{aligned}
 \right .
 \end{equation*}
\begin{equation*}
\iff \left \{
\begin{aligned}
  h_i + \dfrac{1}{2}\norm{\weight_t^{(i,:)}}_1 + \norm{\gain_t^{(i,:)}}_1\varphi_{\max} & \leq U_{\max}  \\
h_i - \dfrac{1}{2}\norm{\weight_t^{(i,:)}}_1 - \norm{\gain_t^{(i,:)}}_1\varphi_{\max} & \geq - U_{\max}  
\end{aligned}
\right .
 \end{equation*}
 \begin{equation*}
\iff \abs{h_i} + \dfrac{1}{2}\norm{\weight_t^{(i,:)}}_1 + \norm{\gain_t^{(i,:)}}_1\varphi_{\max}  \leq U_{\max}.
\end{equation*}
Hence, the constraint \eqref{e:decisionboundsingle} is equivalent to the hard control constraint \eqref{e:controlset}. 
\end{pf}
\begin{pf}[Proof of Theorem \ref{th:DF_rep_stbl}]
\revised{
We consider the first $\reachindex$ blocks of the control $\control_{t:\reachindex}$ according to \eqref{e:decision} and set 
$(\offset_t)_{1:\reachindex m} = - \reachab_{\reachindex}^+ \Aortho^{\reachindex} \sat_{r, \zeta}^\infty \bigl(\stortho_t \bigr)$, $(\gain)_{1:\reachindex m} = \zeros$ and $(\weight)_{1:\reachindex m} = \zeros$. The first claim immediately follows from \cite[Theorem 4]{ref:PDQ-15} and the second claim follows from Lemma \ref{t:msbsingle} and \cite[Theorem 4]{ref:PDQ-15}.
}
\end{pf}  
The proof of Theorem \ref{th:DF_seq} is implied by the proof of Theorem \ref{th:DF_rep}, and the proof of Theorem \ref{th:DF_seq_stbl} is implied by the proof of Theorem \ref{th:DF_rep_stbl}. Details are omitted for the sake of brevity. 

                                 
\end{document}